\documentclass[aps,pre,superscriptaddress]{revtex4}
\usepackage{amsfonts,amsmath,amssymb,amsthm,latexsym,epsfig}

\newtheorem*{teo*}{Theorem}

\newtheorem{lemma*}{Lemma}
\newcommand{\pr}[1]{\left( #1\right)}
\newcommand{\es}[1]{\begin{equation}\begin{split}#1\end{split}\end{equation}}
\newcommand{\est}[1]{\begin{equation*}\begin{split}#1\end{split}\end{equation*}}

\begin{document}

\title{Open Mushrooms: Stickiness revisited}
\author{Carl P. Dettmann \footnote{Carl.Dettmann@bristol.ac.uk}}
\affiliation{School of Mathematics, University of Bristol, United Kingdom}
\author{Orestis Georgiou \footnote{maxog@bristol.ac.uk}}
\affiliation{School of Mathematics, University of Bristol, United Kingdom}
\begin{abstract}
We investigate mushroom billiards, a class of dynamical systems with sharply divided phase space. For typical values of the control parameter of the system $\rho$, an infinite number of marginally unstable periodic orbits (MUPOs) exist making the system sticky in the sense that unstable orbits approach regular regions in phase space and thus exhibit regular behaviour for long periods of time. The problem of finding these MUPOs is expressed as the well known problem of finding optimal rational approximations of a real number, subject to some system-specific constraints. By introducing a generalized mushroom and using properties of continued fractions, we describe a zero measure set of control parameter values $\rho\in(0,1)$ for which all MUPOs are destroyed and therefore the system is less sticky. The open mushroom (billiard with a hole) is then considered in order to quantify the stickiness exhibited and exact leading order expressions for the algebraic decay of the survival probability function $P(t)$ are calculated for mushrooms with triangular and rectangular stems.
\end{abstract}
\maketitle

\section{INTRODUCTION}
Billiards \cite{Tabachnikov} are systems in which a particle alternates between motion in a straight line and specular reflections from the walls of its container. Because they demonstrate a broad variety of behaviours (regular, chaotic \cite{NChe} and mixed phase space dynamics) they have been readily used as models in theoretical and experimental physics \cite{Szasz,Nakano08,Montangero09,Nosich08}.
They are widely applicable because their dynamics corresponds to the classical (short wavelength) limit of wave equations for light, sound or quantum particles in a homogeneous cavity.
\begin{figure}[h]
\begin{center}
\fbox{
\includegraphics[scale=0.2]{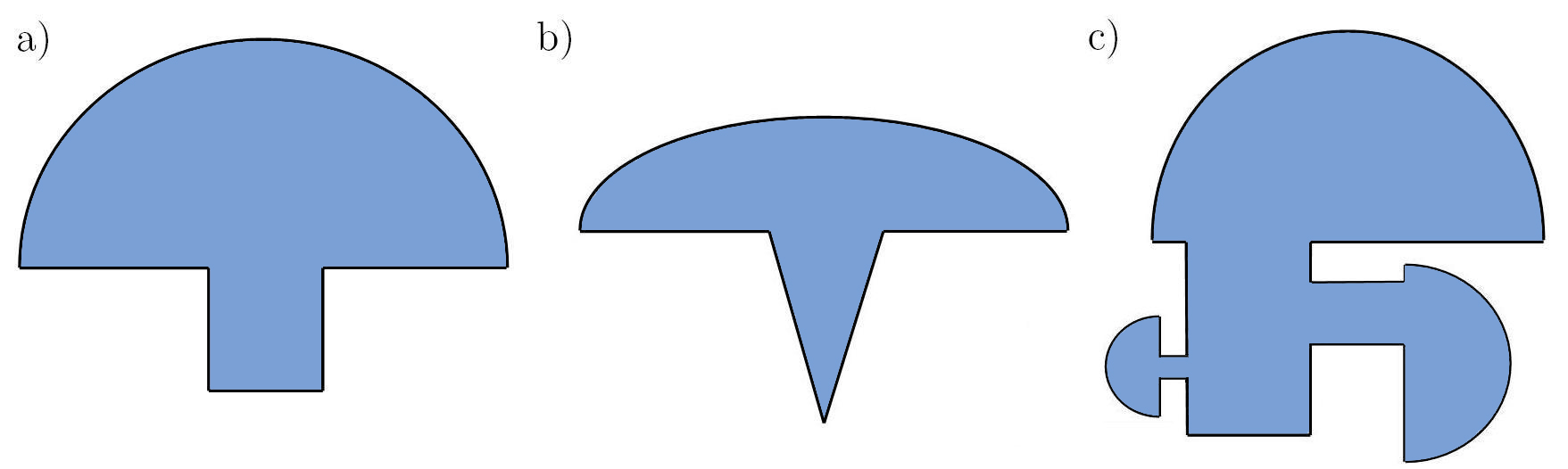}}
\caption{\label{fig:mushrooms} \footnotesize (Color online) a) Simple mushroom, b) Elliptic mushroom with triangular stem, c) `Honey mushroom' with $3$ integrable islands and $1$ ergodic component.}
\end{center}
\end{figure}

The mushroom billiard is constructed by a convex semi-elliptical (including semi-circular) `hat' attached to a `stem' such that their intersection is smaller than the diameter of the hat. Examples are shown in Figure ~\ref{fig:mushrooms}. It is special in that under certain conditions \cite{Bun01,Bun08}, it forms a class of dynamical systems with sharply divided phase space which is easy to visualize and analyse. For example, the phase space of the mushroom shown in Figure ~\ref{fig:mushrooms} a) is composed of a single completely regular (integrable) invariant component and a single connected (topologically transitive) chaotic and ergodic component, in contrast with other generic mixed systems such as the standard map \cite{Ott}, where KAM hierarchical islands form a dense family in the neighborhood of each other. Interestingly, mushrooms can also be designed to have an arbitrary number of integrable and ergodic components (see Figure ~\ref{fig:mushrooms} c)). Hence, mushroom billiards are paradigmatic models for studying the phase space dynamics near the boundary of integrable islands. However, one must note that small perturbations (imperfections) to their boundary may cause the emergence of KAM islands or even complete chaos \cite{Bun08,Seligman08}.

Mushrooms have become increasingly interesting to the quantum chaos community because of their unusually simple, divided classical phase space. This has facilitated for the numerical verification of Percival's conjecture which states that in the semiclassical limit, eigenmodes localize to one or another invariant region of phase space (regular or chaotic), with occurrence in proportion to the respective phase space volumes \cite{Perc73,Barnett07}; recently this has been applied to generalise the boundary term in Weyl's law \cite{Backer10.2}. Similarly, the mechanism of dynamical tunneling between classically isolated phase space regions has also been investigated in the context of mushrooms \cite{Backer08,Backer10} and has been observed in microwave mushrooms \cite{Dietz07}. However, although their classical phase space is sharply divided, generic `simple' mushrooms (Figure ~\ref{fig:mushrooms} a)) have been shown to have an infinite number of different families of marginally unstable periodic orbits (MUPOs), `embedded' in the ergodic component of their phase space \cite{Altman05,Altman06}. These MUPOs live in the mushroom's hat and resemble periodic orbits of the circle. The flow close to these orbits, is reminiscent of that close to KAM islands \cite{Altman08} and thus causes the system to display the phenomenon of `stickiness', where chaotic orbits stick close to regions of stability for long periods of time \cite{Karney83}, the quantum analogue of which is not well understood.

It is worth mentioning that non-sticky mushrooms have been previously constructed using elliptical hats and non-rectangular stems \cite{Bun08}. This is because each focus of the ellipse provides a sharp boundary between rotational and librational orbits, and may be used as the end point of the entrance to the foot. However, in such a case, some care is needed with the stem's length and it's base width, to ensure sufficient defocusing. In addition to this, the size of the opening of the stem must also ensure a bounded number of maximum possible collisions in the hat.

In this paper we focus on classes of mushrooms with circular hats, in which stickiness is due to MUPOs. We express the problem of finding them as the well known problem of finding optimal rational approximations of a number (section II A.). This remarkable connection made with number theory allows us to introduce and characterize a zero measure set of control parameter values, using continued fractions, for which all MUPOs are completely removed (section II B.). This set, not previously discussed in the literature, corresponds to mushrooms with a less sticky hat, the implications of which are yet to be studied classically or quantum mechanically and are hoped to be useful in various applications such as directional emission in dielectric micro-cavities \cite{Ansersen09,Altmann09}. We obtain upper bounds for MUPO-free and finitely sticky irrational mushrooms and also give an explicit example of a MUPO-free mushroom billiard (section II C.).

As discussed in Ref \cite{DettBun07}, placing a small hole on a billiard's boundary allows one to `peep' into the system's dynamics. The smaller the hole, the smaller the observational effect (whether quantum or classical) on the dynamics. We therefore `open' the billiard with a hole and look at the survival probability $P(t)$, given a uniform initial distribution of particles. By considering linear perturbations of MUPOs, we obtain exact expressions for the asymptotic algebraic decay of $P(t)$. This is done for two separate cases, firstly for MUPOs in the semi-circular hat of the mushroom (section III A.) and then for bouncing ball orbits in the case of a rectangular stem (section III B.). The explicit form of these expressions in turn allows us not only to predict but also to calibrate the survival probability function by changing the various geometric parameters involved. Finally, the results are confirmed numerically (section III C.) and then discussed briefly, including higher dimensional mushrooms and other implications of this work (section IV).

\section{STICKINESS IN CLOSED MUSHROOMS}

\subsection{MUSHROOMS WITH MUPOS}

\begin{figure}[h]
\begin{center}
\fbox{
\includegraphics[scale=0.2]{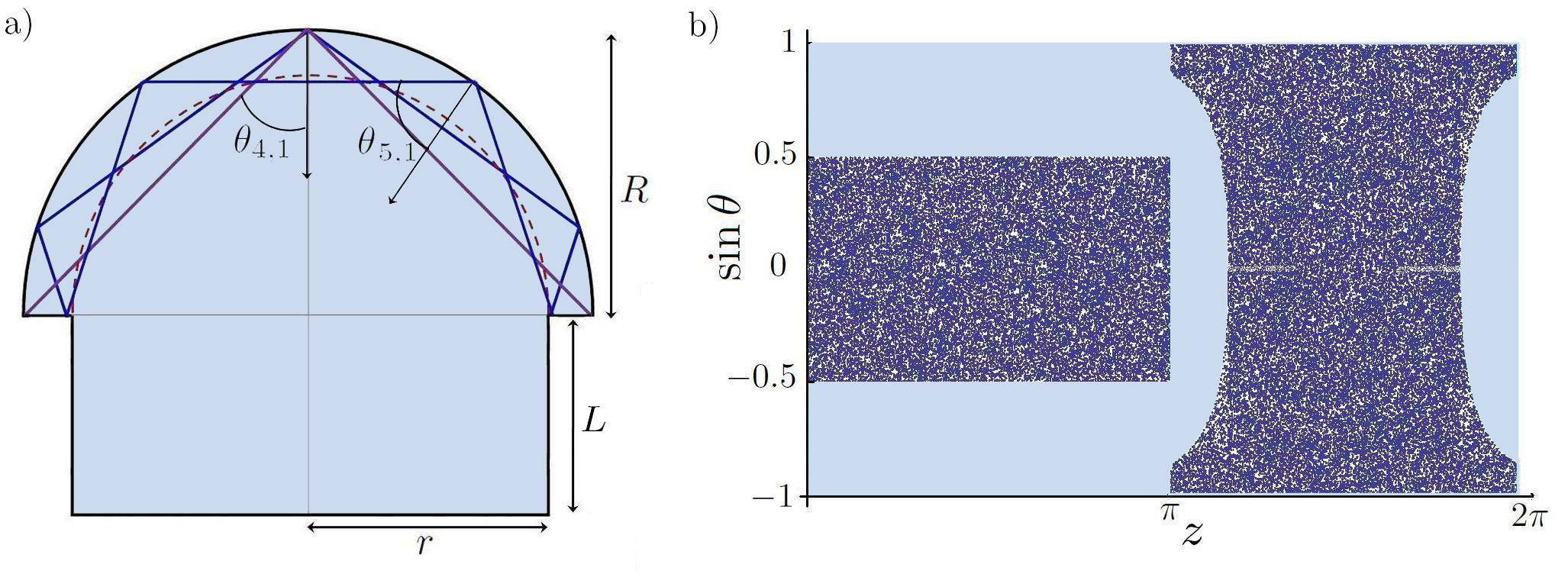}}
\caption{\label{fig:mushroom} \footnotesize (Color online) \emph{Left}: Any orbit intersecting the dashed red semicircle of radius $r$ is unstable while any orbit not intersecting it is stable. MUPOs are the periodic orbits which intersect the semicircle while not entering the mushroom's stem. The MUPOs $(s,j)=(4,1)$ and $(s,j)=(5,1)$ are shown. \emph{Right}: Phase space plot of the ergodic component for a simple mushroom using Birkhoff coordinates ($z,\sin\theta$), with $\frac{r}{R}=0.5$ and $L=0.5$ computed from a random chaotic orbit followed for $10^{5}$ collisions. Here, $z\in[0,\pi R+ 2(R+L))$ is the arc length parametrization along the billiard's boundary, increasing from zero at from the right-most point of the mushroom, in an anticlockwise fashion and $\theta\in\pr{-\frac{\pi}{2},\frac{\pi}{2}}$ is the angle of incidence at each collision.}
\end{center}
\end{figure}

Stickiness \cite{Karney83} is due to points or regions of phase space having vanishingly small local Lyapunov\cite{Yao93} exponents, typically found in the vicinity of KAM elliptic islands or in the very close neighborhood of marginally unstable periodic orbits (MUPOs). Sticky orbits exhibit long periods of quasi-regular behaviour and are thus also associated with the phenomenon of intermittency \cite{Schuster}. A well known and well studied family of MUPOs are the bouncing ball orbits present in the chaotic Bunimovich stadium billiard \cite{Bu74}.
They form a zero measure family of period two orbits \cite{Berry78}, trapped forever between the stadium's parallel walls. Even though MUPOs do not affect the overall ergodicity of the system, they govern long time statistical properties of the system, such as the Poincar\'{e} recurrence times distribution $\mathcal{Q}(t)\sim t^{-2}$ \cite{Altman09}, the rate of mixing (the rate of the decay of correlations) $\mathcal{C}(t)\sim t^{-1}$ \cite{Markarian04,Chernov05,Chernov08} and the long time survival probability $P(t)\sim t^{-1}$ \cite{OreDet09} given a carefully positioned `hole'. Furthermore, the exponents of these power-laws appear to be a universal fingerprint of nonuniform hyperbolicity and stickiness, at least for one and two dimensional Hamiltonian systems with sharply divided phase space \cite{Altman06}. The MUPOs in the mushroom's hat and in the annular billiard were extensively studied by Altmann in his PhD thesis \cite{AltmanPhD} and more briefly in \cite{Altman05,Altman06,Altman08} and are common in many billiards with circular arcs, the significance of which has only recently been realised in the context of directional emission in dielectric micro-cavities \cite{Ansersen09,Altmann09}.

MUPOs in the mushroom are best understood when introduced geometrically. The dashed red semicircle of radius $r$ in Figure ~\ref{fig:mushroom}a) corresponds to the border between the ergodic and regular component of the mushroom's phase space (see Figure ~\ref{fig:mushroom}b). Any orbit intersecting this semicircle is unstable and lies in the ergodic component in phase space while any orbit not intersecting it exhibits regular motion and remains forever in the mushroom's hat. However, as shown in the figure, one can find periodic orbits in the mushroom's hat which do intersect it and therefore are unstable. A compact way of describing them is given by:
\begin{equation}
\alpha_{s,j}= \cos\frac{j \pi}{s}\leq\frac{r}{R} < \frac{\cos\frac{j \pi}{s}}{\cos\frac{\pi}{\lambda s}}=\beta_{s,j},
\end{equation}
where
\begin{equation}
s \geq 3,
\qquad
1 \leq j \leq \begin{cases}
\frac{s}{2}-1, &   \text{if $s$ is even},\\
\frac{s-1}{2}, &  \text{if $s$ is odd}, \end{cases}
\qquad
\lambda =\begin{cases}
1, & \text{if $s$ is even},\\
2, &  \text{if $s$ is odd}.\end{cases}
\end{equation}
In Eq(1), $r$ and $R$ are as defined in Figure ~\ref{fig:mushroom}. The coprime integers $s$ and $j$ describe periodic orbits of the circle billiard with angles of incidence $\theta_{s,j}=\frac{\pi}{2}-\frac{j \pi}{s}$. More specifically, $s$ is the period and $j$ the rotation number of the orbit. $R\alpha_{s,j}$ is the shortest distance from the periodic orbit  $(s,j)$ to the origin. $R\beta_{s,j}$ is half the longest straight line passing through the origin which intersects the unfolded (along the hat's base) periodic orbit $(s,j)$ at equal distances on either side. Hence, (1) guarantees that $(s,j)$ is a MUPO and can be oriented in such a way as not to enter the stem while still intersecting the dashed semicircle. Let $\mathcal{S}_{\rho}$ denote the set of periodic orbits which are marginally unstable for a given $\rho=\frac{r}{R}$.

It can easily be seen that a small perturbation $\eta$ with respect to the incidence angle $\theta_{s,j}$ of a MUPO will cause the orbit to precess in the opposite direction and eventually enter the stem and feel the chaotic effect of the defocusing mechanism \cite{Wojt86}. However, since the precessing angular velocity is proportional to the perturbation strength $\eta$ which may be arbitrarily small, the orbit will behave in a quasi-periodic fashion and entry into the stem may take an unbounded amount of time. Hence the term `stickiness', meaning that orbits in the immediate vicinity of MUPOs stick close to the regular component of phase space for long periods of time.
Note however that although these periodic orbits are dynamically marginally unstable, they are not structurally robust against parameter perturbations of $\rho$.

It is clear from (1) that the intervals $(\alpha_{s,j},\beta_{s,j})$ are shrinking quadratically with increasing $s$. We see this by rearranging (1) into:
\begin{equation}
\frac{j}{s}\geq \frac{1}{\pi}\arccos\rho > \frac{1}{\pi}\arccos\left(\frac{\cos\frac{j \pi}{s}}{\cos\frac{\pi}{\lambda s}}\right),
\end{equation}
expanding for large $s$
\begin{equation}
\frac{j}{s}\geq  \vartheta^{*} > \frac{j}{s} - \left(\frac{\pi \cot \frac{j \pi}{s}}{2 }\right)\frac{1}{\lambda^2 s^2} + \mathcal{O}\left(\frac{1}{s^4}\right),
\end{equation}
where we have set $\vartheta^{*}=\frac{1}{\pi}\arccos\rho$ and rearranging once more to get
\begin{equation}
0\leq \frac{j}{s}- \vartheta^{*} < \left(\frac{\pi \cot \frac{j \pi}{s}}{2 \lambda^2}\right)\frac{1}{s^2},
\end{equation}
where we have neglected the positive terms of order $\sim s^{-4}$, thus possibly losing some of the MUPOs; we give explicit bounds on this term in the next section.
In this way the problem of finding the elements of $\mathcal{S}_{\rho}$ is expressed as the well known number theoretic problem of finding rational approximations $\frac{j}{s}$ of $\vartheta^{*}\in(0,\frac{1}{2})$. However in this case we have a couple of complications: the approximations are one-sided and the tolerance depends both on the numerical value of $\frac{j}{s}$ and the parity of $s$ through $\lambda$ in (2).

This interesting connection made here allows one to apply well known results from number theory to the present dynamical system and infer useful dynamical properties about it. Altmann \cite{AltmanPhD} proved that for almost all $\rho\in(0,1)$ there exist infinitely many MUPOs and hence the corresponding mushroom exhibits stickiness and all inherent dynamical properties mentioned above. In the following section we remove the parity dependence in the context of a more general mushroom model. This in turn allows us to use properties of continued fractions more effectively to derive a sufficient condition so that (1) has no solutions and hence destroy all MUPOs in the hat of the mushroom.

\subsection{MUSHROOMS WITHOUT MUPOS}
\subsubsection{\textbf{Generalized Mushroom}}
\begin{figure}[h]
\begin{center}
\fbox{
\includegraphics[scale=0.3]{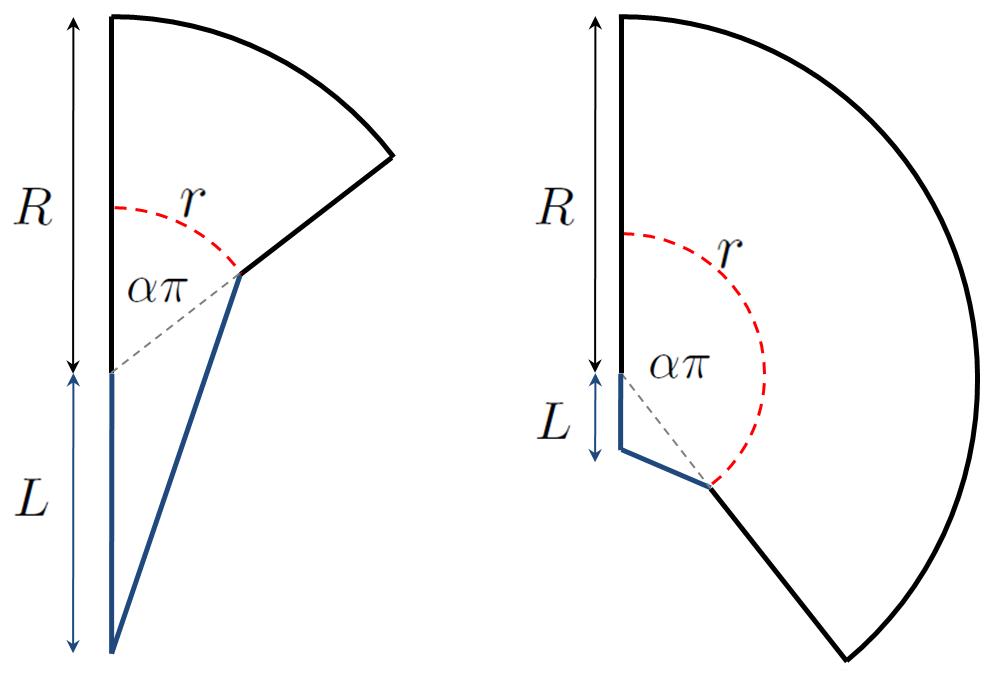}}
\caption{\label{fig:mushroom2} \footnotesize (Color online) Generalized mushroom billiards with variable hat sizes and triangular stems.}
\end{center}
\end{figure}

We have seen that MUPOs in the mushroom billiard, which give it its sticky dynamical character, are directly related to number theory through (5). In this section we propose a generalization of the mushroom billiard studied in section II A., which will allow us to efficiently use properties of continued fractions without having to worry about the parity of MUPOs. The main result here will be to prove the existence of a zero measure set of $\rho=\frac{r}{R}$ values for which the mushroom's hat is MUPO-free. Furthermore, we shall obtain a sufficient condition which explicitly describes a subset of this set.

Consider the `elementary cell' obtained by slicing the mushroom along its vertical axis of symmetry. Then the period of the corresponding $(s,j)$ orbit is $s \lambda /2$. Similarly, since we are currently only interested in collisions with the curved segment of the billiard, we introduce the parameter $\alpha\in(0,1)$ which allows the mushroom to have circular hats of variable size.
This billiard, shown in Figure ~\ref{fig:mushroom2}, is in a class of billiards considered in \cite{Dietz06}, and was shown to have a sharply divided phase space in \cite{Bun08} as long as $L>0$. The boundary between the two components is given by the dashed arc of radius $r\in(0,R)$. Notice that because the stem is triangular, there are no bouncing ball orbits present.

Periodic orbits in the hat of the proposed mushroom will now have incidence angles with the curved boundary given by $|\theta_{q,p}|=\frac{\pi}{2}-\frac{\alpha p \pi}{q}$ for some coprime $p$ and $q$, and equation (1) becomes
\begin{equation}
\cos\frac{\alpha p \pi}{q}\leq \frac{r}{R} < \frac{\cos\frac{\alpha p \pi}{q}}{\cos\frac{\alpha \pi}{q}},
\end{equation}
Notice that there is no longer a parity dependent $\lambda$. Similarly (5) becomes
\begin{equation}
0\leq \frac{p}{q}- \frac{\vartheta^{*}}{\alpha} < \frac{\alpha \pi \cot \frac{\alpha p \pi}{q}}{2 q^{2}}+ \frac{R_{2}(q)}{q^{2}},
\end{equation}
where $\vartheta^{*}=\frac{1}{\pi}\arccos\rho$, and $R_{2}(q)$ is the remainder term obtained from the Taylor expansion for large $q$. In the following we bound the argument of the cotangent by $\pi\vartheta^{*}$ and bound $R_{2}$, so that for $q\geq Q$ we have
\begin{equation}
0\leq \frac{p}{q}- \frac{\vartheta^{*}}{\alpha} < \frac{\alpha \pi \cot \frac{\alpha p \pi}{q}}{2 q^{2}}+ \frac{R_{2}(q)}{q^{2}}<\bigg[ \overbrace{\bigg(\frac{\alpha \pi \rho}{2 \sqrt{1 - \rho^2}}\bigg)+\hat{R}_{2}(q,Q)}^{K(q,Q)}\bigg]\frac{1}{q^2},
\end{equation}
where $R_{2}(q)$ is bounded by
\es{
\hat{R}_{2}(q,Q)= \frac{\alpha^{2} \pi^{2}}{\cos^{2}\pr{\frac{\alpha \pi}{Q}}q^{2}}\Bigg[&\left(\tan^{2}\pr{\frac{\alpha \pi}{Q}}+\frac{4}{3}\right)\frac{\rho}{\sqrt{1-\rho^{2}}} + \frac{\rho^{3}}{2 (1-\rho^{2})^{\frac{3}{2}}}
\\ &+ \pr{1+\frac{\alpha^{2}\pi^{2}}{\cos(\frac{\alpha\pi}{Q})^{2}Q^{2}}}\frac{\rho^3}{2 \cos(\frac{\alpha \pi}{Q})^{2}\pr{1-\rho^2(1+\frac{\alpha^2\pi^2}{\cos(\frac{\alpha \pi}{Q})^{2}{q^{2}}})^2}^\frac{3}{2}}\Bigg],
}
as obtained in Appendix \textit{A.} Here, $Q$ is a fixed number up to which (8) must be checked numerically. It must be greater than $\mathrm{max}\pr{\alpha \pi, \frac{\alpha \pi}{\cos 1}\sqrt{\frac{\rho}{1-\rho}}}$, following from Appendix A.

\subsubsection{\textbf{MUPO-free Condition}}

We now turn to some number theory and introduce some basic concepts. It is well known that the best rational approximations of a real number $\xi$ are obtained through its continued fraction representation \cite{Szusz}
\begin{equation}
\xi= a_{0}+ \frac{1}{a_{1}+\frac{1}{a_{2}+ \ldots}}=[a_{0};a_{1},a_{2}, \ldots],
\end{equation}
where the quantities $a_{0},a_{1},a_{2},\ldots$ are called `partial quotients' and are usually taken to be positive integers. Irrational numbers have an infinite continued fraction representation while rationals have finite. The $n^{th}$ truncation of a continued fraction representation gives the $n^{th}$ `convergent' $\frac{A_{n}}{B_{n}}$ of $\xi$. Hence irrational numbers have an infinite number of convergents while rationals finite. Convergents are `best approximations' to $\xi$, meaning that there is no other fraction with denominator smaller than $B_{n}$ which approximates $\xi$ better.

For the mushroom, we would like to find values of $\frac{\vartheta^{*}}{\alpha}=[a_{0};a_{1},a_{2}, \ldots]$ for which
\begin{equation}
0\leq \frac{p}{q}- \frac{\vartheta^{*}}{\alpha} < \frac{K(Q,Q)}{q^{2}},
\end{equation}
has no solutions since this would also imply no solutions to (6). Solutions to (11), if any, are only given by the convergents of $\frac{\vartheta^{*}}{\alpha}$ if $0<K(Q,Q)\leq\frac{1}{2}$ \cite{Khinchin97}. Obviously, if $\frac{\vartheta^{*}}{\alpha}$ is rational then there is only a finite number of solutions to (11). However if $\frac{\vartheta^{*}}{\alpha}$ is irrational the answer is not so simple. We focus on the convergents $\frac{p}{q}=\frac{A_{n}}{B_{n}}$ of $\frac{\vartheta^{*}}{\alpha}$ and express it in terms of them such that
\begin{equation}
\frac{\vartheta^{*}}{\alpha}= \frac{\zeta_{n+1}A_{n} +A_{n-1}}{\zeta_{n+1} B_{n} +B_{n-1}},
\end{equation}
where $\zeta_{n}=[a_{n};a_{n+1},a_{n+2}, \ldots]$ is the $n^{th}$ `complete quotient' of $\frac{\vartheta^{*}}{\alpha}$. Hence
\begin{equation}
\frac{A_{n}}{B_{n}}- \frac{\vartheta^{*}}{\alpha}= \frac{A_{n}B_{n-1}-A_{n-1}B_{n}}{\left(\zeta_{n+1}+\frac{B_{n-1}}{B_{n}}\right)B_{n}^{2}} = \frac{(-1)^{n-1}}{\left(\zeta_{n+1}+\frac{B_{n-1}}{B_{n}}\right)B_{n}^{2}}.
\end{equation}
It is easy to see that if $n$ is even, then $\frac{A_{n}}{B_{n}}- \frac{\vartheta^{*}}{\alpha}<0$. Therefore, equation (11) will not have any solutions if
\begin{equation}
K(Q,Q)<\frac{1}{\zeta_{n+1}+\frac{B_{n-1}}{B_{n}}},
\end{equation}
for all odd $n$. Since $a_{n+1}<\zeta_{n+1}<a_{n+1}+1$ and $\frac{B_{n-1}}{B_{n}}< 1$, it follows that
\begin{equation}
K(Q,Q)<\frac{1}{\wp+2},
\end{equation}
where $\wp=max(a_{2 n})$, is a sufficient condition for (11) and therefore (6) not to have any solutions. The condition is never satisfied if $a_{2 n}$ is unbounded.

The set of numbers with bounded even partial quotients as derived above has zero measure \cite{Khinchin97} and has Hausdorff dimension one as $\wp$ is unbounded as $\rho\rightarrow 0$ \cite{Hensley92}. As shown by Altmann \cite{AltmanPhD}, a generic mushroom will have infinitely many MUPOs and so it will be infinitely sticky. However, we have shown here that there are infinitely many values of $\frac{\vartheta^{*}}{\alpha}$, and therefore $\rho$, for which MUPOs in the hat are finite or completely removed. Since the smallest possible value of $\wp$ is one, for the original mushroom with $\alpha= \frac{1}{2}$, if $\rho< \left(\left(\frac{3\pi\alpha}{2}\right)^{2}+1\right)^{-\frac{1}{2}}\approx 0.390683$, (15) gives a sufficient condition for (11) not to have any solutions and therefore describes a mushroom with no MUPOs in its hat.

\subsection{MUSHROOMS WITH LARGE STEMS}
\subsubsection{\textbf{MUPO-free Example}}
For larger values, $\frac{1}{2}<K(Q,Q)\leq 1$, solutions to (11), if any, are given by the convergents $\frac{A_{n}}{B_{n}}$ and also by the so called `intermediate' convergents of the form $\frac{c A_{n+1}+A{n}}{c B_{n+1}+B_{n}}$ \cite{Khinchin97}, where $c$ is an integer such that $1\leq c< a_{n+2}$.
There are however values of $\vartheta^{*}$ with $K(Q,Q)>\frac{1}{2}$ satisfying (14) such that the corresponding mushrooms will have no MUPOs. An example of such a mushroom is $\rho=\cos\left(\frac{5+\sqrt{2}}{23}\pi\right)\approx0.64013$ which has $K(q,95)< 0.6549$ and $2\vartheta^{*}=[0; 1, 1, 3, \{1, 4\}]$ (where we have numerically checked the absence of MUPOs up to $q=95$). Here, the odd convergents of $2\vartheta^*$ satisfy $0\leq \frac{A_{n}}{B_{n}}- \frac{\vartheta^{*}}{\alpha}= \frac{K_{n}}{B_{n}^{2}}$ where $K_{n}=\left(\zeta_{n+1}+\frac{B_{n-1}}{B_{n}}\right)^{-1}$ for odd $n\geq3$, where $\zeta_{n+1}=[1;4,\{1,4\}]=\frac{1}{2}(1+\sqrt{2})$. It is an easy exercise to show that for all odd $n\geq3$
\es{B_{n}&=\frac{1}{4} \pr{\alpha_{-}\lambda_+^\frac{n}{2} - \alpha_{+}\lambda_-^\frac{n}{2}},\\
B_{n-1}&= \frac{1}{8} \pr{\beta_{-} \lambda_+^\frac{n}{2} + \beta_{+} \lambda_-^\frac{n}{2}},
}
where $\lambda_\pm=3\pm2 \sqrt{2}$, $\alpha_{\pm}=12\pm7\sqrt{2}$ and $\beta_{\pm}=\pm26+19\sqrt{2}$ are all positive numbers. Hence $\frac{B_{n-1}}{B_{n}}= \frac{1}{2}\pr{\frac{\beta_{-}+\beta_{+}\pr{\frac{\lambda_{-}}{\lambda_{+}}}^{\frac{n}{2}}}{\alpha_{-}-\alpha_{+}\pr{\frac{\lambda_{-}}{\lambda_{+}}}^{\frac{n}{2}}}}$ is strictly decreasing with $n$ and therefore $K_{n}$ is bounded by
\es{K(q,95)<K_{5}\leq K_{n}<\frac{1}{\sqrt{2}},}
for all odd $n\geq5$, where $K_{5}\approx0.706$. Similarly for the intermediate convergents of $2\vartheta^{*}$ we have that
\es{\frac{c A_{n+1}+A_{n}}{c B_{n+1}+B_{n}}- 2\vartheta^{*}= \frac{c A_{n+1}+A_{n}}{c B_{n+1}+B_{n}} - \frac{\zeta_{n+2} A_{n+1}+ A_{n}}{\zeta_{n+2}B_{n+1}+B_{n}},
}
which for odd $n\geq5$ simplifies to
\es{ \frac{2+2\sqrt{2}-c}{\pr{c B_{n+1}+B_{n}}\pr{B_{n+1}(2+2\sqrt{2})+B_{n}}}\equiv \frac{\bar{K}_{n}(c)}{(c B_{n+1}+B_{n})^{2}},
}
since $\zeta_{n+2}=[4;1,\{4,1\}]=2+2\sqrt{2}$. Hence, using (16) and a similar argument as above $\bar{K}_{n}(c)= \frac{4+4c-c^{2}}{4 \sqrt{2}}- \frac{(2+5\sqrt{2})(c-2-2\sqrt{2})^{2}}{8(5-\sqrt{2})} \pr{\frac{\lambda_{-}}{\lambda_{+}}}^{\frac{n}{2}}$ is bounded by
\es{ K(q,95) < \bar{K}_{5}(1) \leq \bar{K}_{n}(c)<\frac{4+4c-c^{2}}{4 \sqrt{2}},
}
for $c=1,2,3$ and odd $n\geq5$ where $\bar{K}_{5}(1)\approx1.237$. Therefore, $\rho=\cos\left(\frac{5+\sqrt{2}}{23}\pi\right)$ describes a mushroom with no MUPOs in its chaotic region.

\subsubsection{\textbf{Supremum of MUPO-free Values}}

From the example above we can now use similar arguments to establish that MUPO-free values of $\rho$ exist up to $\frac{1}{\sqrt{2}}$. In other words $\mathrm{sup}\big(\rho\in(0,1):\mathcal{S}_{\rho}=\emptyset\big)=\frac{1}{\sqrt{2}}$. To see this, let $\hat{K}(Q,Q)$ denote the value of $K(Q,Q)$ at $\rho=\frac{1}{\sqrt{2}}$. Then from equations (8) and (9) $K(Q,Q)< \hat{K}(Q,Q)=\frac{\pi}{4}+\frac{7\pi^{2}}{12 Q^{2}}+\mathcal{O}(Q^{-4})$ for $0<\rho<\frac{1}{\sqrt{2}}$. Now consider for $m\in \mathbb{Z}^{+}$ large
\es{2\vartheta^{*}=[0;1,\{1,m\}]= \frac{m+2+\sqrt{m^{2}+4m}}{4m-4}=\frac{1}{2}+\frac{1}{4 m}+\mathcal{O}\pr{\frac{1}{m^{2}}},
}
so that $\rho=\cos \pi\vartheta^{*}= \frac{1}{\sqrt{2}}-\frac{\pi}{8\sqrt{2}m}+\mathcal{O}(m^{-2})$. We first look at the odd convergents of $2\vartheta^{*}$ as in (13)
\es{0<\frac{A_{n}}{B_{n}}-2\vartheta^{*} \equiv \frac{K_{n}}{B_{n}^{2}},}
where $K_{n}=\pr{\zeta_{n+1}+\frac{B_{}n-1}{B_{n}}}^{-1}$ and $\zeta_{n+1}=[1;m,\{1,m\}]=\frac{1}{2}+\sqrt{\frac{1}{4}+\frac{1}{m}}$. Via a similar manipulation as in (16) we obtain that for odd $n\geq3$
\es{\frac{B_{n-1}}{B_{n}}=\frac{(2-\lambda_{-})+(\lambda_{+}-2)\Big(\frac{\lambda_{-}}{\lambda_{+}}\Big)^{\frac{n-1}{2}}}{(1+2m-\lambda_{-})-(1+2x-\lambda_{+})\Big(\frac{\lambda_{-}}{\lambda_{+}}\Big)^{\frac{n-1}{2}}}
,}
where $\lambda_{\pm}=\frac{1}{2}\big(2+m\pm\sqrt{4m^{2}+m^{3}}\big)$. Thus $K_{n}\geq K_{1}= \frac{2 m}{3m +\sqrt{m(4+m)}}=\frac{1}{2}-\frac{1}{4 m}+\mathcal{O}(m^{-2})$ converges exponentially to $\frac{m}{\sqrt{m(4+m)}}=1-\frac{2}{m}+\frac{6}{m^{2}}+ \mathcal{O}(m^{-3})$ with $n$ and therefore $K_{n}>\hat{K}(q,B_{n})\approx\frac{\pi}{4}$ for large enough $m$ and $n$. Similarly, when looking at the intermediate convergents of $2\vartheta^{*}$ as in (18) and (19) such that
\es{\bar{K}_{n}(c)=\frac{(\zeta_{n+2}-c)\Big(c+\frac{B_{n}}{B_{n+1}}\Big)}{\Big(\zeta_{n+2}+\frac{B_{n}}{B_{n+1}}\Big)},
}
where $c=1,2,\ldots(m-1)$, $\zeta_{n+2}=[m;1,\{m,1\}]=\frac{m}{2}+\sqrt{\frac{m^{2}}{4}+m}$ and $\frac{B_{n}}{B_{n+1}}$ can be obtained from (23), we find that $\bar{K}_{n}(1)\leq \bar{K}_{n}(c)$ and $\bar{K}_{n+1}(c)<\bar{K}_{n}(c)$. Hence, since $\bar{K}_{n}(1)$ converges exponentially to $\frac{1-2m}{m^{2}-m\sqrt{m(4+m)}-1}=1-\frac{2}{m^{2}}+\mathcal{O}(m^{-3})$ with $n$, then $\bar{K}_{n}(c)\geq\bar{K}_{n}(1)> \frac{1-2m}{m^{2}-m\sqrt{m(4+m)}-1} >\hat{K}(q,c B_{n+1}+B_{n}) \approx\frac{\pi}{4}$ for large enough $m$ and $n$, thus verifying our claim above for the supremum of MUPO-free mushrooms.

\subsubsection{\textbf{Supremum of Finitely Sticky Irrational Values}}

Another claim we can make is that $\mathrm{sup}\pr{\rho: \# \mathcal{S}_{\rho}<\infty, \frac{1}{\pi}\arccos\rho \not\in \mathbb{Q}}=\frac{4}{\sqrt{16+\pi^{2}}}\approx 0.7864$. To see this we take the leading order term of $K(Q,Q)$ as $Q\rightarrow\infty$ and equate it to one, so that $\frac{\pi \rho}{4 \sqrt{1-\rho^{2}}}=1$.
Now since $\rho=\cos \pi \vartheta^{*}$, then $2\vartheta^{*}=\frac{2}{\pi}\arccos\frac{4}{\sqrt{16+\pi^{2}}}=[0;2,2,1,3,1,1,1,\ldots]<[0;2,2,1,3,1,2,1,\{1,m\}]=2\tilde{\vartheta}^{*}$.
It follows that the corresponding $\tilde{K}(Q,Q)<1$ in the limit $Q\rightarrow\infty$. Now since we may augment the tail of the continued fraction expansion of $2\vartheta^{*}$ as done above by the transformations $a_{\nu}\rightarrow a_{\nu}+1$ and $\zeta_{\nu+2}\rightarrow [1;m,\{1,m\}]$ for any even $\nu$, then $K_{n}\equiv B_{n}^{2}\pr{\frac{A_{n}}{B_{n}}-2\tilde{\vartheta}^{*}}$ will converge exponentially to some function $f(m)= 1-\frac{k}{m}+ R_{1}(m)$ with $n$ for some constant $k$ and $R_{1}(m)= \frac{f''(\xi)}{2 m^{2}}$ for some $0<\xi<m$. Therefore, as $\nu\rightarrow \infty$, $(\tilde{\vartheta}^{*}-\vartheta^{*})\rightarrow 0^{+}$ and $\tilde{K}(Q,Q)\rightarrow 1$. However we may always choose $m$ and $n$ big enough such that $K_{n}>\tilde{K}(B_{n},B_{n})$. A similar statement can be made for the intermediate convergents of $2\tilde{\vartheta}^{*}$. For values of $\rho>\frac{4}{\sqrt{16-\pi^{2}}}$, $K(Q,Q)>1$ and therefore all convergents of $2\vartheta^{*}\not\in\mathbb{Q}$ are solutions of (11) \cite{Szusz} hence describing mushrooms with infinitely many MUPOs.

In the next section we investigate the stickiness in the hat of the mushroom in the context of escape through a small hole placed on the stem of the mushroom.

\section{ESCAPE FROM THE MUSHROOM}

In the previous sections we have investigated the dynamics of the mushroom billiard and specifically focused at the chaotic region of phase space, close to the regular island. We have seen how MUPOs come into existence, how they affect the dynamics of orbits in their immediate vicinity but also how they can be removed. In this section we shall address the problem of escape through a small hole placed on the stem of the mushroom. We shall derive exact expressions, to leading order, for the survival probability for two specific cases of the mushroom.

The uniform (Liouville) distribution projected onto the billiard boundary has the form $(2|\partial Q|)^{-1}\mathrm{d}z$ $\mathrm{d}\sin \theta $, where $|\partial Q|$ is the perimeter of the billiard while $z\in(0,|\partial Q|]$ is the length parametrization round the billiard boundary and $\theta\in(-\frac{\pi}{2},\frac{\pi}{2})$ is the angle of incidence with it. This is the most natural choice for an initial distribution of particles. Given such a distribution, the probability $P(t)$ that a particle survives (\textit{i.e.} does not escape through $k$ small holes) in a fully chaotic billiard up to time $t$ decays exponentially $\sim e^{-\gamma t}$, with $\gamma\approx \frac{\sum_{i=1}^{k} \epsilon_{i}}{\langle\tau\rangle |\partial Q|}$ \cite{DettBun07} to leading order, where $\langle\tau\rangle$ is the mean free path, $\epsilon_{i}$ is the length of each hole and $|Q|$ the area of the billiard. $|\partial Q|$ is easy to calculate and hence for the remainder of the paper we leave it in general form to allow for different shaped stems.

For billiards with mixed phase space such as the mushroom, $P(t)$ for long enough times is expected to decay as
\begin{equation}
P(t)\approx \mathcal{A} + \mathcal{B}\left( e^{-\bar{\gamma} t} + \frac{\mathcal{C}}{t}\right),
\end{equation}
where we have neglected terms of order $t^{-2}$. $\mathcal{A}$ is the measure of the integrable island given by:
\begin{equation}
\mathcal{A}=4 (2|\partial Q|)^{-1}\left[R\sqrt{1-\rho^{2}}-\rho R \arccos \rho+ \frac{\pi}{2}R(1-\rho)\right],
\end{equation}
and $\mathcal{B}$ is its complement ($\mathcal{B}= 1-\mathcal{A}$).
In (25) we have assumed that the holes are placed well in the ergodic component of phase space and therefore
\begin{equation}
\bar{\gamma}\approx\frac{\sum_{i=1}^{k} \epsilon_{i}}{\langle\bar{\tau}\rangle \mathcal{B}|\partial Q|}
\end{equation}
while the mean free path in the ergodic component is now
\begin{equation}
\langle\bar{\tau}\rangle= \frac{c_{\nu}}{c_{\mu}} = \frac{\pi\left(|Q_{s}| + R^{2}\arcsin \rho +\rho R^{2} \sqrt{1-\rho^{2}}\right)}{\mathcal{B} |\partial Q|},
\end{equation}
where $|Q_{s}|$ is the area of the mushroom's stem while $c_{\nu}$ and $c_{\mu}$ are the invariant probability measures of the ergodic component for the billiard flow and map respectively.

Algebraic decays, of the form $\frac{\mathcal{C}}{t}$ in (25), originate from the stickiness exhibited, which is itself due to the nonuniform hyperbolicity of intermittent systems as discussed in the previous sections. It is a geometrical description of the constant $\mathcal{C}$ that we seek here. In the case of the stadium billiard for example, near-bouncing ball orbits were studied and such a constant was successfully calculated in Ref \cite{OreDet09}. In the following two subsections we attempt to do the same, first for the MUPOs living the mushroom's hat and then for near-bouncing ball orbits present in mushrooms with rectangular stems.

\subsection{STICKY HAT}
\begin{figure}[h]
\begin{center}
\fbox{
\includegraphics[scale=0.3]{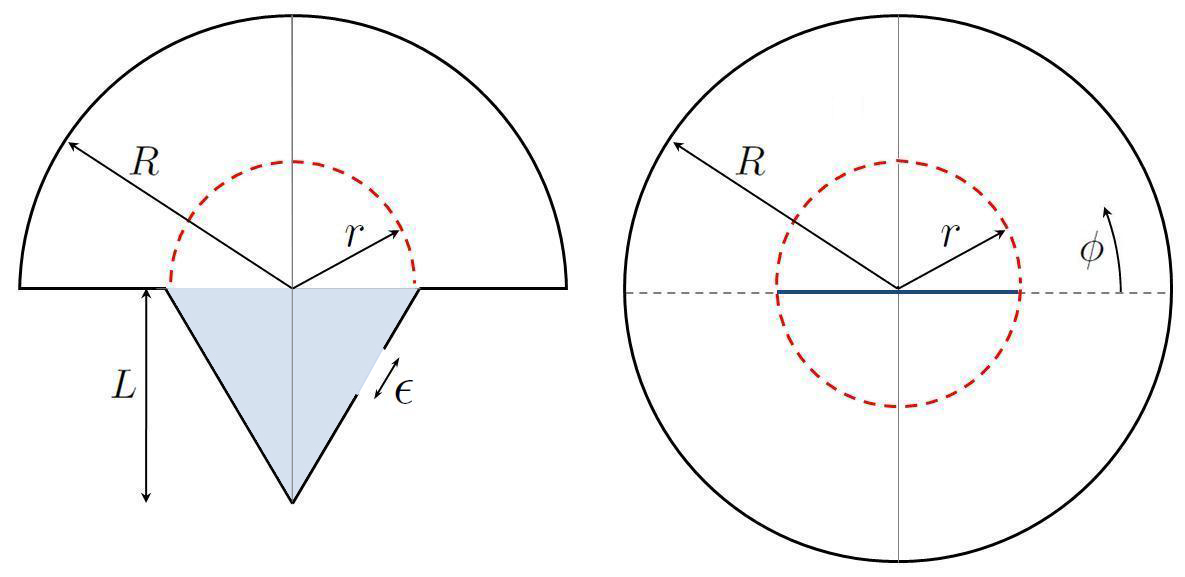}}
\caption{\label{fig:mushroom3} \footnotesize (Color online) \emph{Left}: Mushroom with triangular stem. \emph{Right}: Image reconstruction trick at the base of the mushrooms hat. Orbits entering the lower semicircle through the thick blue horizontal line of length $2r$ are assumed to escape through the hole $\epsilon$ soon thereafter.}
\end{center}
\end{figure}

We consider a mushroom with a central triangular stem and circular hat as shown in Figure ~\ref{fig:mushroom3}, hence removing any bouncing ball orbits between parallel walls. The only source of stickiness is hence found in the immediate vicinity of the MUPOs in the hat and therefore the algebraic decay $\frac{\mathcal{C}}{t}$ is equal to the measure (relative volume occupied in phase space) of the set of quasi-periodic initial conditions which do not enter the stem until a time $t$. This approximation is well justified since any orbit entering the stem has a very small probability of re-entering a `sticky' mode before suffering escape according to the exponential decay of (25).

We use the \textit{image reconstruction trick} \cite{Bun01} and neglect collisions with the base of the mushroom's hat. Hence the dynamics in the hat remains unchanged while a horizontal slit of length $2r$ centered at the origin corresponds to the stem's opening. We parametrize the now circular boundary by the angle $\phi$, where $\phi\in(0,2 \pi)$ increases anticlockwise as shown in the right panel of Figure ~\ref{fig:mushroom3}. Now, it is easy to see that each initial condition $(\phi, \theta_{s,j})$ is a MUPO if the collision coordinate $\phi$ satisfies:
\begin{equation}
\phi\in\bigcup_{k=0}^{\lambda s-1}\left(\phi_{1}(\theta_{s,j},k),\phi_{2}(\theta_{s,j},k)\right),
\end{equation}
with
\begin{align}
\phi_{1}(\theta_{s,j},k)&= \theta_{s,j} +\frac{\pi}{\lambda}+ \arccos\left(\rho^{-1}\sin\theta_{s,j}\right)+ (k-1) \frac{2 \pi}{\lambda s},\\
\phi_{2}(\theta_{s,j},k)&= \theta_{s,j} +\frac{\pi}{\lambda}- \arccos\left(\rho^{-1}\sin\theta_{s,j}\right)+ k \frac{2 \pi}{\lambda s},
\end{align}
where $\rho=\frac{r}{R}$ and the angles $\phi_{i}$ are taken modulo $2\pi$.
Each MUPO then defines a dashed, horizontal line in the $\phi-\theta$ plane (the phase space), and each dashed line has length $\phi_{2}- \phi_{1}= \frac{2\pi}{\lambda s}- 2 \arccos\left(\rho^{-1}\sin \theta_{s,j}\right)$. Notice that $\phi_{1}$ and $\phi_{2}$ are not defined if $\sin\theta > \rho$.

\begin{figure}[h]
\begin{center}
\fbox{
\includegraphics[scale=0.4]{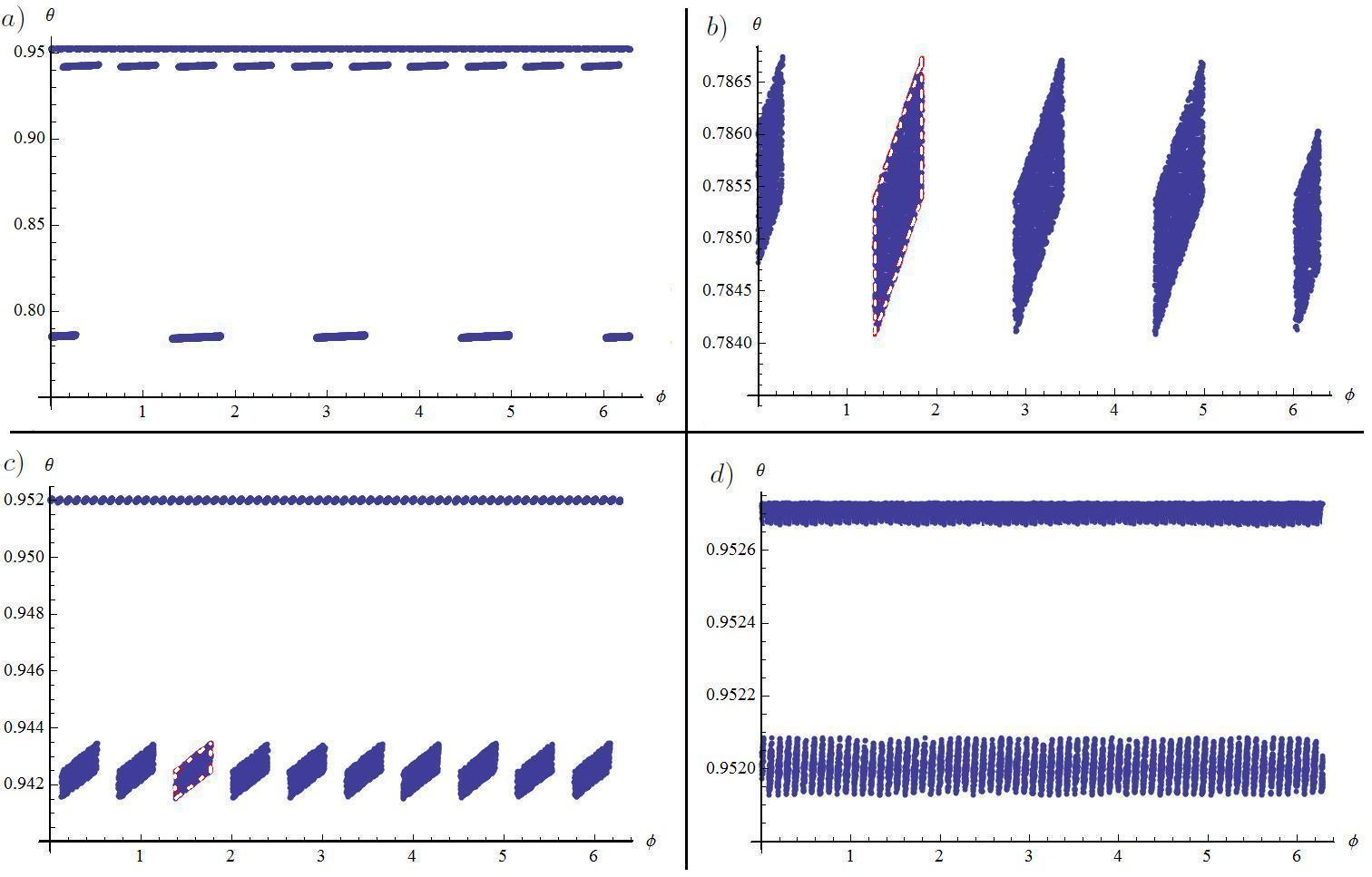}}
\caption{\label{fig:phase} \footnotesize (Color online) Phase space plots of the Circle billiard with a vertical slit hole in the center (see also Figure ~\ref{fig:mushroom3}). The dynamics are same as for the hat of the mushroom billiard in the sense that the ICs shown in blue will not escape through the slit and hence not enter the stem of the mushroom for at least $N$ collisions with the boundary. In all plots above we have used $\rho= 0.815$ and $N=200$. Plots $b)$ $c)$ and $d)$ are magnifications of $a)$, showing in more detail
the MUPOs (4,1), (5,1), (66,13) and their surrounding sticky orbits. It is clear that MUPOs accumulate closer to the boundary of the
integrable island at $\theta=\arcsin\rho$. The white dashed lines in $b)$ and $c)$ are the analytic prediction given by (29-32).)}
\end{center}
\end{figure}

To help visualize how the long surviving initial conditions near the above described MUPOs populate the phase space, we turn to some computer simulations. Initial conditions near the integrable island's boundary are chosen randomly so that $\phi\in(0,2\pi)$ and $\theta\in(0,\arcsin\frac{r}{R})$. The ones that survive for at least $N$ collisions with the boundary are shown in the top left panel of Figure ~\ref{fig:phase} for parameters $N=200$ and $\rho= 0.815$. We notice that for the selected value of $\rho$, the most dominant MUPO is the square with $(s,j)=(4,1)$ (see also Figure ~\ref{fig:phase}b)). In Figure ~\ref{fig:phase}c) one can identify the pentagon orbit $(s,j)=(5,1)$ which like all odd $s$-orbits has twice its period ($\lambda s=10$) of surviving intervals along the horizontal line $\theta_{5,1}$. Further magnification into the phase space reveals the $(s,j)=(66,13)$ orbit (see also Figure ~\ref{fig:phase}d)) and then an accumulation of higher order orbits closer to the island's boundary at $\arcsin\rho$. The next MUPO is $(s,j)=(920,181)$.

We introduce a small perturbation $\eta \ll 1$ in the angle $\theta_{s,j}$ of each MUPO and expand (30-31) to leading order:
\begin{align}
\phi_{1}(\theta_{s,j}+\eta,k)= \phi_{1}(\theta_{s,j},k)+ \left(1-\frac{ \cos\theta_{s,j}}{\sqrt{\rho^{2}-\sin^{2}\theta_{s,j}}}\right)\eta+\mathcal{O}(\eta^{2}),\\
\phi_{2}(\theta_{s,j}+\eta,k)= \phi_{2}(\theta_{s,j},k)+\left(1+\frac{ \cos\theta_{s,j}}{\sqrt{\rho^{2}-\sin^{2}\theta_{s,j}}}\right)\eta+\mathcal{O}(\eta^{2}).
\end{align}
We also impose a time constraint such that the perturbed MUPO will survive up to time $t$ by requiring that
\begin{align}
\phi \geq  \phi_{1}(\theta_{s,j}+\eta,k) + 2 \eta N,\\
\phi \leq  \phi_{2}(\theta_{s,j}+\eta,k) + 2 \eta N.
\end{align}
where $N =\frac{t}{2 R \cos(\theta_{s,j}+\eta)}$ is the number of collisions in time $t$. Expanding (34-35) to leading order defines in total $4$ lines which form a quadrilateral in phase space with area $\Delta_{s,j}$ which can be integrated with respect to the invariant measure $(2|\partial Q| \mathcal{B})^{-1}\mathrm{d}\phi $ $ \mathrm{d}\sin \theta $ to give:
\begin{equation}
\Delta_{s,j}= \frac{8 R \cos^{2}\theta_{s,j}\left(\pi - s \lambda \arccos\left(\rho\sin\theta_{s,j}\right)\right)^2}{2 s^{2}\lambda^{2} |\partial Q| \mathcal{B} t} +\mathcal{O}\left(\frac{1}{t^2}\right),
\end{equation}
to leading order in $t$. There are $2 \lambda s$ such quadrilaterals due to $\theta$-symmetry, however only half of the total area for each MUPO lies in $\phi\in(0,\pi)$, which corresponds to the actual mushroom's hat. As for the initial conditions on the straight segments of the hat, since the billiard map is measure preserving, only $2\lambda j$ quadrilaterals are mapped onto them. Hence overall we obtain:
\begin{equation}
\frac{\mathcal{C}}{t}= \sum_{(s,j)\in\mathcal{S}_{\rho}} \lambda(s + 2 j) (\Delta_{s,j}-\delta_{s,j}) + \mathcal{O}\left(\frac{1}{t^{2}}\right),
\end{equation}
where $\mathcal{S}_{\rho}$ was defined in section II A. and
\begin{equation}
\delta_{s,j}= \begin{cases}
\Delta_{s,j}/2, &   \text{ if } \cos\frac{j \pi}{s}=\rho ,\\
0, &  \text{otherwise}.\end{cases}
\end{equation}
accounts for the possibility that a MUPO is situated exactly on the border of the chaotic region and therefore can only be perturbed from one side. The sum in (37) converges since the elements of $\mathcal{S}_{\rho}$, if any, are distributed with a bounded density with respect to $\ln s$. Also, notice that $\mathcal{C}$ does not depend on the size or position of the hole on the stem. Numerical simulations of the survival probability function are performed and discussed in section III C.

\subsection{STICKY STEM}

\begin{figure}[h]
\begin{center}
\fbox{
\includegraphics[scale=0.33]{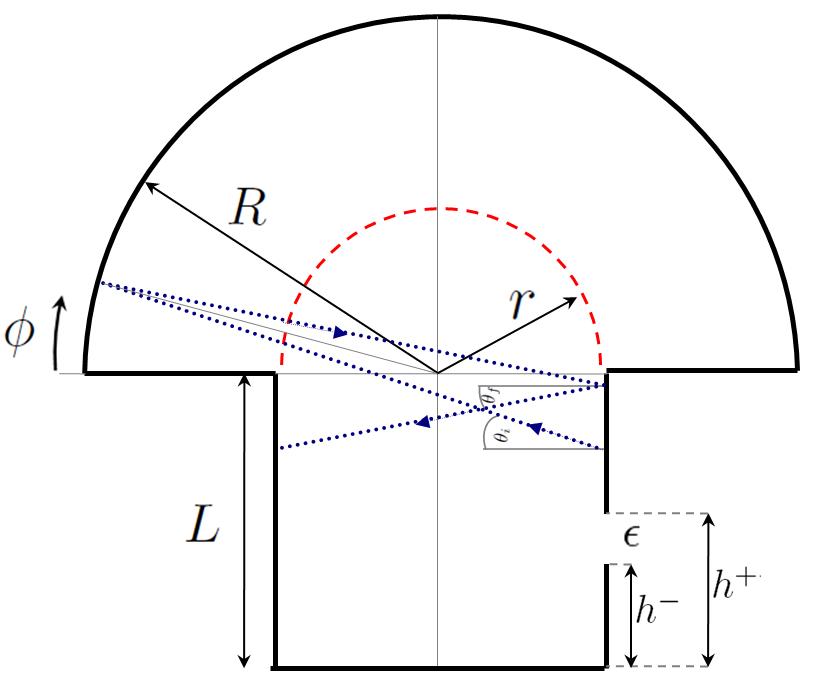}}
\caption{\label{fig:mushroom4} \footnotesize (Color online) Mushroom with rectangular stem and a hole of size $\epsilon$ on one of its parallel walls. A near-bouncing ball orbit experiencing a nonlinear collision process in the mushroom's hat is shown.}
\end{center}
\end{figure}

In the previous section we derived an expression to leading order for the asymptotic behaviour of $P(t)$ (see equations (25-29) and (36-38)) for a sticky mushroom with a triangular stem. Here we investigate the stickiness introduced by the bouncing ball orbits present in mushrooms with rectangular stems of length $L$ and a hole of size $\epsilon$ on one of the two parallel segments as shown in Figure ~\ref{fig:mushroom4}. A method for calculating the contribution of these orbits to $P(t)$ was devised and explained in detail in \cite{OreDet09}. Here, we follow this method and obtain an exact expression to leading order for the survival probability of the mushroom billiard. In doing so we discover an interesting discontinuous dependence of $P(t)$ on $\rho=\frac{r}{R}$ and also show that in the limit $\rho\rightarrow 1$ the expression for $P(t)$ reduces to the one obtained in \cite{OreDet09} for the stadium billiard.

We first split the billiard's boundary $\partial Q$ into four, non-overlapping, connected segments: $\partial Q_{s}^{b}, \partial Q_{s}^{w}, \partial Q_{h}^{b}$ and $\partial Q_{h}^{c}$, referring to the stem's base, the stem's parallel walls, the hat's base and the hat's curved segment respectively.
We parametrize the right parallel wall of $\partial Q_{s}^{w}$ by $s\in(0,L)$ such that the interval $(h^{-},h^{+})$ defines the hole of size $\epsilon$ as shown in Figure ~\ref{fig:mushroom4}. It is now easy to see that initial conditions (ICs) $(x,\theta)$ with $0<x<h^{-}$ and $0<\theta<\arctan\frac{\epsilon}{4\rho R}$ cannot jump over the hole and therefore do not interact with the mushroom's hat. Such orbits behave in a completely regular manner and therefore can be integrated directly to give
\begin{equation}
\frac{2}{ 2|\partial Q| \mathcal{B}} \Big(\int_{0}^{\arcsin\left(h^{-}/t\right)}\int_{t \sin\theta }^{h^{-}}\cos\theta \mathrm{ d}s \mathrm{ d}\theta + \int_{0}^{\arcsin\left(2 h^{-}/t\right)}\int_{t \sin\theta }^{2 h^{-}}\cos\theta \mathrm{ d}s \mathrm{ d}\theta \Big)=\frac{(2h^{-})^{2}+(h^{-})^{2}}{2|\partial Q| \mathcal{B}t},
\end{equation}
where we have neglected terms of order $\sim t^{-2}$ and multiplied by $2$ due to the horizontal symmetry of the billiard.
Similarly, ICs with $h^{+}<x<L$ and $0>\theta>-\arctan\frac{\epsilon}{4\rho R}$ give
\begin{equation}
\frac{2}{ 2|\partial Q|\mathcal{B}} \int_{0}^{\arcsin\left((L-h^{+})/t\right)}\int_{t \sin\theta }^{L-h^{+}}\cos\theta \mathrm{ d}s \mathrm{ d}\theta= \frac{(L-h^{+})^{2}}{2|\partial Q|\mathcal{B}t}.
\end{equation}
ICs from $\partial Q_{s}^{b}$ have contributions of order $\sim t^{-2}$ to $P(t)$ and therefore are ignored.

As expected, the survival probability at long times is proportional to the square of the available length on either side of the hole. For the remainder of this section we consider ICs $(x_{i},\theta_{i})$ such that $h^{+}<x_{i}<L$ and $0< \theta_{i} \ll 1$, and investigate how they contribute to $P(t)$. We let $n$ denote the number of collisions a particle experiences from straight to straight segment before entering the hat of the mushroom, and define $d_{1}= L-(x_{i} + 2 r n \tan\theta_{i})>0$ as the distance from the edge of the straight to the point of the last straight wall collision. We can see that $n=\Big\lfloor \frac{L-x_{i}}{2 r \tan\theta_{i}}\Big\rfloor$, where $\lfloor \cdot \rfloor$ $\lceil\cdot\rceil$ are the floor and ceiling functions respectively. Note that $0< d_{1} < 2 r \tan\theta_{i}$. Once a particle enters the hat of the mushroom it is advantageous to switch to coordinates suitable for the circle billiard map given by $(\phi,\psi)\rightarrow(\phi+\pi-2\psi,\psi)$ such that $\phi$ is the angular collision coordinate and increases from zero in an anticlockwise fashion as shown in Figure ~\ref{fig:mushroom4}, while $\psi\in(-\frac\pi2,\frac\pi2)$ is the angle of reflection. Note that $\phi$ is different from what was used in section III A. Also, $\psi$ is used instead of $\theta$ here to distinguish between collisions on the curved segment of the billiard boundary ($\partial Q_{h}^{c}$) and collisions elsewhere. Once in the hat, we neglect collisions with the vertical base $\partial Q_{h}^{b}$, by using the image reconstruction trick as before.
We find that the particle entering the hat will first collide with $\partial Q_{h}^{c}$ at
\begin{equation}
\phi=-\frac{d_{1}}{R}+\left(1+\rho\right)\theta_{i}>0,
\end{equation}
and its angle will be
\begin{equation}
\psi=-\frac{d_{1}}{R}+\rho\theta_{i}.
\end{equation}
Let $\theta_{f}$ be the final angle obtained when the orbit re-enters the stem of the mushroom after experiencing a reflection
process (a series of $k\in\mathbb{Z}^{+}$ collisions with $\partial Q_{h}^{c}$) in the hat. We thus find that
\begin{equation}
\theta_{f}= \frac{2 k d_{1}}{R}-\left(2 k \rho + 1\right)\theta_{i}.
\end{equation}
By carefully investigating the reflection process we find that $k$ is actually restricted to only three possible scenarios such
that $k$ can either be equal to $1$, $\left\lceil\frac{R}{r}\right\rceil$ or $\left\lceil\frac{R}{r}\right\rceil+1$, depending on
the ICs ($x_{i},\theta_{i}$), which agree with the so called `magic numbers' from Ref \cite{Seligman08}. This can be seen if one looks at the least number of iterations of the
circle billiard map before the orbit described by (41) and (42) intersects the horizontal slit hole:
\begin{equation}
k= inf\left\{j\in\mathbb{Z}^{+}:\left|\frac{\psi}{\left(2j-1\right)\psi+\phi}\right|<\rho\right\}.
\end{equation}

\begin{figure}[h]
\begin{center}
\fbox{
\includegraphics[scale=0.3]{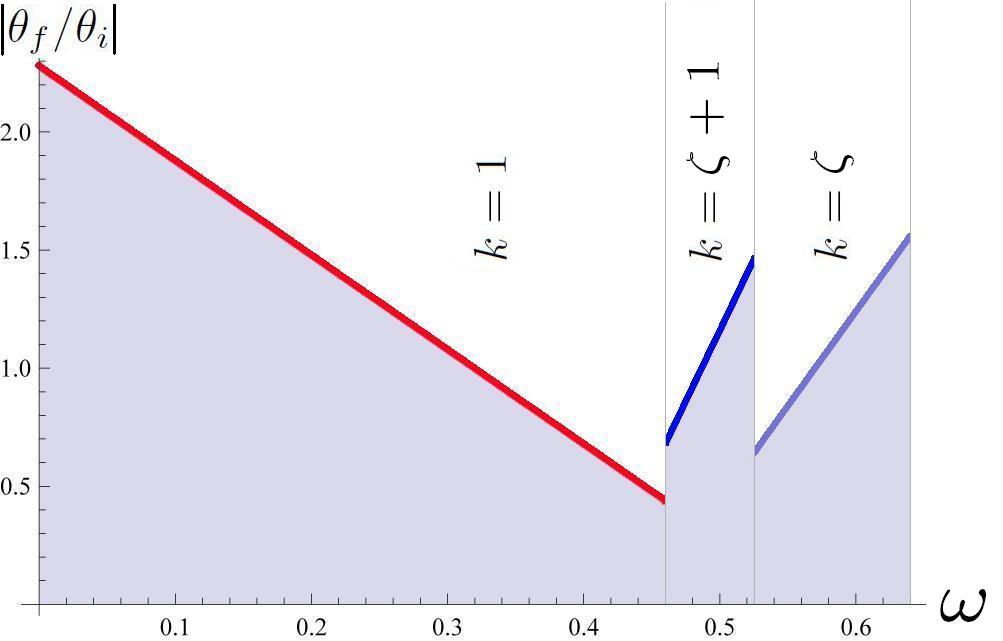}}
\caption{\label{fig:process} \footnotesize (Color online) Reflection process in the mushroom's hat described in (43) using the substitution $\omega=\frac{d_{1}}{2 R \theta_{i}}\in(0,\rho)$ with $\rho=\cos\frac{(5+\sqrt{2})\pi}{23}$ and $d_{1}= 0.01$. The Red, Blue and light Blue lines correspond to $k= 1 , \zeta+1$ and $\zeta$ reflection process respectively as inscribed in the figure.}
\end{center}
\end{figure}

In equation (45) below we have substituted the possible values of $k$ into (43) and also calculated the values of $\theta_{i}$ for
which each collision scenario corresponds to:
\begin{equation}
\theta_{f}=\begin{cases}
\frac{2 d_{1}}{R}-\left(2 \rho + 1\right)\theta_{i}<0, \hspace{33 mm}\frac{(2\rho +1)d_{1}}{2\rho(\rho+1)R} <\theta_{i}, \hspace{31 mm} k=1 \text{ collision }\\
\\
\frac{2 \zeta d_{1}}{R}-\left(2 \zeta \rho + 1\right)\theta_{i}>0, \hspace{37 mm}\frac{d_{1}}{2\rho R} < \theta_{i}< \frac{(2\zeta \rho -1)d_{1}}{2\zeta \rho^{2}R}, \hspace{12.5 mm} k=\zeta \text{ collisions }\\
\\
\frac{2 \left(\zeta+1\right) d_{1}}{R}-\left(2 \left(\zeta+1\right) \rho + 1\right)\theta_{i}>0, \hspace{13.5 mm}\frac{(2\zeta \rho -1)d_{1}}{2\zeta \rho^{2}R}<\theta_{i}<\frac{(2\rho +1)d_{1}}{2\rho (\rho+1)R}, \hspace{13 mm} k=(\zeta+1) \text{ collisions }.\end{cases}
\end{equation}
where we have set $\zeta=\left\lceil\frac{R}{r}\right\rceil$. Note that $d_{1}$ is a function of both $x_{i}$ and $\theta_{i}$. The first inequality on $\theta_{i}$ ($k=1$ collision) seems
to suggest that $\theta_{i}$ is unbounded, however this is not the case. This can be seen in an example situation plotted in Figure ~\ref{fig:process} where we have made the substitution $\omega=\frac{d_{1}}{2 R \theta_{i}}\in(0,\rho)$.
Notice that if $\rho^{-1}$ is an integer, then the $\zeta$-collision processes in (45) is no longer attainable
and we only have two possible collision scenarios. It is interesting to note that if $\rho=1$, equation (45) reduces to equations
(10-11) of \cite{OreDet09} which refer to the stadium billiard's reflection process with the curved segment.

We now formulate the time of escape for ICs $(x_{i},\theta_{i})$:
\begin{equation}
t(x_{i},\theta_{i},k)\approx\frac{L-x_{i}}{\theta_{i}}+ \frac{L-h^{+}}{\left|\theta_{f}\right|}+ 2 R\left(\rho + k+1\right),
\end{equation}
where we have taken small angle approximations, and substitute the values of $\theta_{f}$ and $k$ for each collision scenario to get three equations for the time to escape.
Each one of these equations describes conic sections since they are quadratic in both $x_{i}$ and $\theta_{i}$ variables. Rearranging
to make $\theta_{i}$ the subject, we obtain three hyperbolae in the $x_{i}-\theta_{i}$ plane, describing the ICs that escape exactly at
large times $t$. It is important to know the domain of each hyperbola. This can be obtained by substituting for the $d_{1}$ variables
into the inequalities of (45), and then rearranging for $\theta_{i}$. These inequalities are given below for the corresponding collision scenarios:
\begin{align}
&k=\zeta \text{ collisions }, \hspace{50 mm}\frac{L - x_{i}}{2 \rho(1 + n)R } <\theta_{i}< \frac{(2 \zeta \rho - 1) (L - x_{i})}{2 \rho \zeta \rho R(1 + 2 n) - 2 \rho n R},\nonumber \\
&k=\left(\zeta+1\right) \text{ collisions }, \hspace{22.5 mm}\frac{(2 \zeta \rho - 1) (L - x_{i})}{2 \rho \zeta \rho R(1 + 2 n) - 2 \rho n R}<\theta_{i}<\frac{(2 \rho + 1) (L - x_{i})}{2 \rho R (\rho + 2 n \rho + 1 + n)},\nonumber\\
&k=1 \text{ collision }, \hspace{35.3 mm}\frac{(2 \rho + 1) (L - x_{i})}{2 \rho R (\rho + 2 n \rho + 1 + n )} <\theta_{i}.
\end{align}

\begin{figure}[h]
\begin{center}
\fbox{
\includegraphics[scale=0.32]{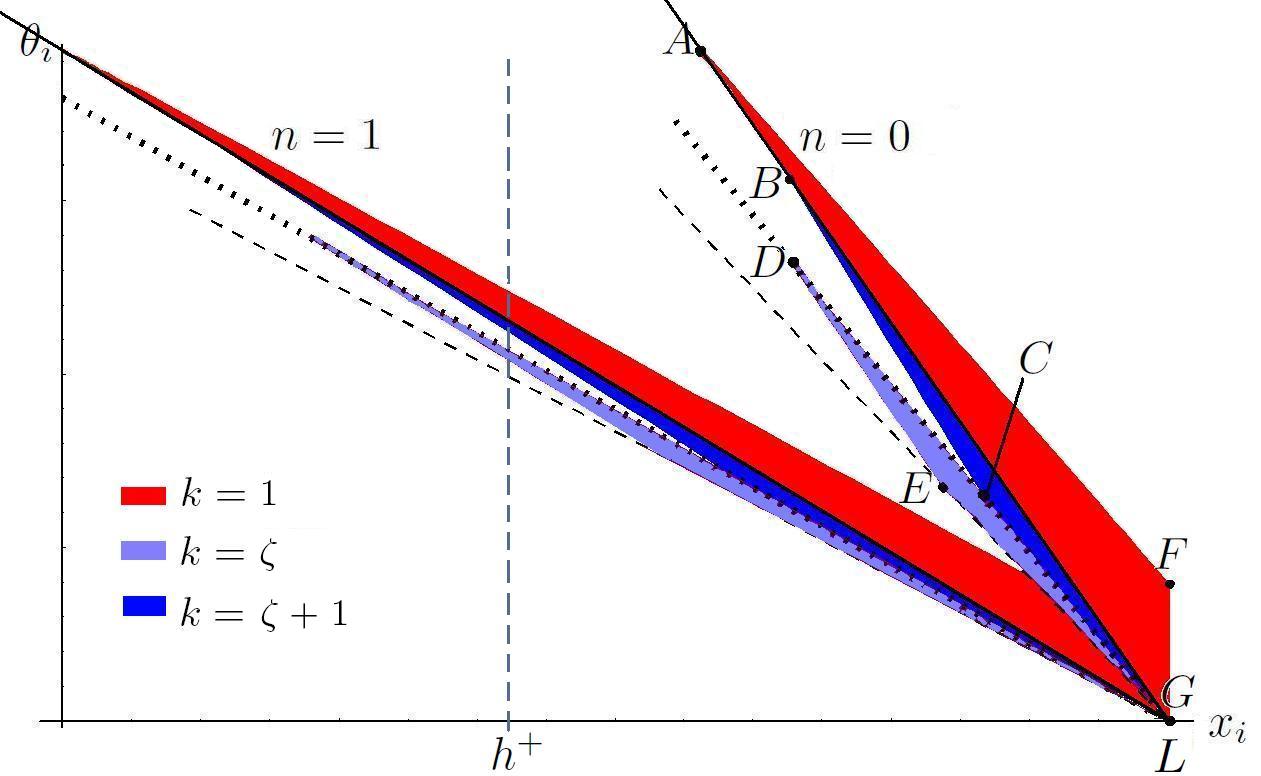}}
\caption{\label{fig:hyperbolas} \footnotesize (Color online) Area enclosed by equations (46) subject to (47) and equations (45) in the $x_{i}-\theta_{i}$ phase space for $n=0$ and $1$, using $\rho=0.6$, $L=1$ and $t=50$. The colours used are in consistence with the ones in Figure ~\ref{fig:process}. The dotted, dashed and solid black straight lines come from the inequalities in equation (45). The corners $A-G$ are defined in Appendix \textit{B.} and are each highlighted by a black dot for $n=0$. The dashed vertical line at $x_{i}=h^{+}$ shows how the hole truncates the area of interest. The area defined for all $n$, corresponds to the ICs that survive at least until time $t=50$.}
\end{center}
\end{figure}

For $n=0, 1, \ldots$ and for $t$ large, we plot the three hyperbolae from (46) subject to (47) and the three straight lines from (45) onto the $x_{i}-\theta_{i}$ plane (see Figure ~\ref{fig:hyperbolas}). These define an area in phase space which corresponds to the ICs that survive at least until time $t$ for fixed $n$. The various colors indicate the type of reflection process $k$ the ICs experience in consistence with the ones in Figure ~\ref{fig:process}.
Notice that as the number of collisions $n$ with the straight segments increases, the area of interest tilts and stretches in a
non-overlapping fashion. To obtain the contribution to $P(t)$ of these long surviving ICs, we must integrate each non-overlapping area and sum them all up. Note that the invariant measure will be assumed to be $\mathrm{d} \mu= (2 |\partial Q| \mathcal{M})^{-1} \mathrm{d} \theta_{i} \mathrm{d} x_{i}$ here since $\theta_{i}$ is small and thus $\text{ d}\sin\theta_{i}\approx\text{ d}\theta_{i}$.

The corners of each enclosed area $A-G$, as shown in Figure ~\ref{fig:hyperbolas}, for each value of $n$ are given in Appendix \textit{B}. There are various issues which one needs to consider in order to obtain correct asymptotic expressions for the areas. Firstly, one needs to approximate all the hyperbolae by straight lines. This is done by joining the corners $A-G$ and thus forming an irregular polygon. For example, the hyperbola between $A$ and $F$, which comes from $t(x_{i},\theta_{i},1)$, is approximated by a straight line joining $A$ and $F$. Similarly, for the hyperbola joining $B$ and $C$, which comes from
$t(x_{i},\theta_{i},\zeta)$, and for the hyperbola joining $D$ and $E$ which comes from $t(x_{i},\theta_{i},\zeta+1)$. The remaining edges are already straight lines and thus need no approximating. As shown in \cite{OreDet09}, the error in these approximations is $\mathcal{O}\left(t^{-2}\right)$ and hence meets our required asymptotic accuracy.

Another issue to be dealt with is the position of the hole which restricts the irregular polygons in $x_{i}\in(h^{+},L)$. This forces a deformation by truncating each polygon from the left each time one of its corners surpasses the hole's position as seen for example in Figure ~\ref{fig:hyperbolas}. This is due to the
tilting effect caused as $n$ is increased. Following \cite{OreDet09} again, we expect $7$ different sums since there are $6$ corners ($A-F$), each of which will intersect the hole at $h^{+}$ at different values of $n$.
We thus solve for $n$ and find that the leftmost corner $A_{x_{i}}=h^{+}$  when $n=n_{A}=\left\lfloor \frac{\frac{t}{R}\left(1+2\rho\right)-2-20\rho-12\rho^{2}-4\rho^{3}}{4\rho+12\rho^{2}+8\rho^{3}}\right\rfloor$. Similar expressions have been obtained for all other corners ($B-F$) and are given in Appendix \textit{C}.
\begin{figure}[h]
\begin{center}
\fbox{
\includegraphics[scale=0.25]{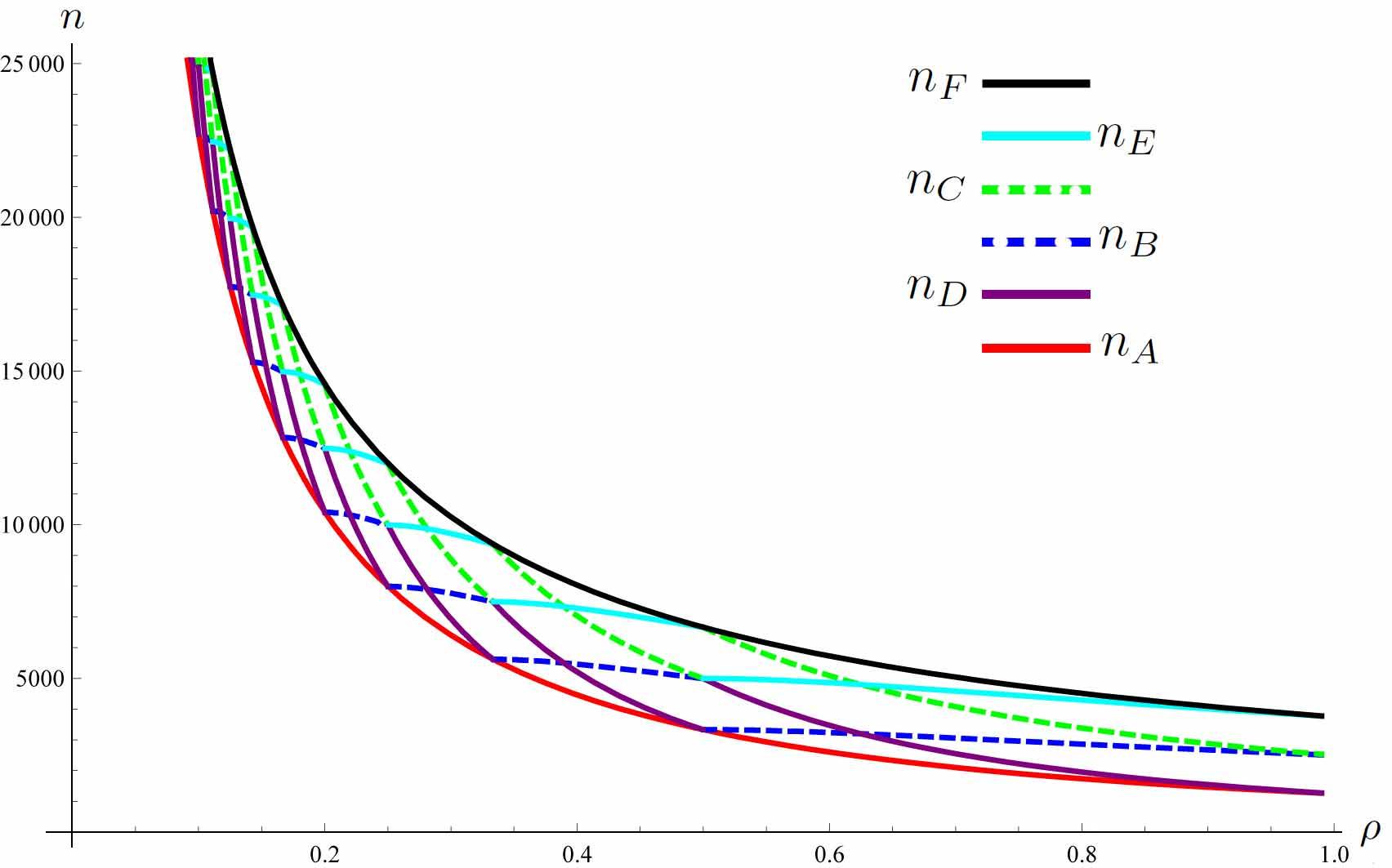}}
\caption{\label{fig:wire} \footnotesize (Color online) $n_{A}-n_{F}$ are defined in the text above and given in Appendix \textit{D}. The figure shows how they vary discontinuously as a function of $\rho=\frac{r}{R}\in(0,1)$ for $t=10^{4}$.}
\end{center}
\end{figure}
Interestingly, we find that the order in which the corners $A-F$ coincide with the hole's position depends on the system's control parameter $\rho=\frac{r}{R}$. Their order alternates between $n_{A}<n_{B}\leq n_{D}<n_{E}\leq n_{C}<n_{F}$ and $n_{A}<n_{D}\leq n_{B}<n_{C}\leq n_{E}<n_{F}$ for $\rho\in(0,1)$, which is shown in Figure ~\ref{fig:wire}. This is due to the discontinuity introduced by the ceiling function in $\zeta$ for $n_{B}-n_{E}$, hence the lower and upper
bounds of the $7$ sums will depend on the above order, and so will their arguments. Altogether we write:
\begin{equation}
\sum_{n=0}^{n_{A}} \hat{P}_{1} + \sum_{n=n_{A}+1}^{n_{B}} \hat{P}_{2}+ \sum_{n=n_{B}+1}^{n_{D}} \hat{P}_{3}+\sum_{n=n_{D}+1}^{n_{E}} \hat{P}_{4}+\sum_{n=n_{E}+1}^{n_{C}} \hat{P}_{5}+\sum_{n=n_{C}+1}^{n_{F}} \hat{P}_{6}+\sum_{n=n_{F}+1}^{\infty} \hat{P}_{7} ,
\end{equation}
\begin{equation}
\sum_{n=0}^{n_{A}} \tilde{P}_{1} + \sum_{n=n_{A}+1}^{n_{D}} \tilde{P}_{2}+ \sum_{n=n_{D}+1}^{n_{B}} \tilde{P}_{3}+\sum_{n=n_{B}+1}^{n_{C}} \tilde{P}_{4}+\sum_{n=n_{C}+1}^{n_{E}} \tilde{P}_{5}+\sum_{n=n_{E}+1}^{n_{F}} \tilde{P}_{6}+\sum_{n=n_{F}+1}^{\infty} \tilde{P}_{7},
\end{equation}
where $\hat{\cdot}$ and $\tilde{\cdot}$ are used to distinguish between the two orderings described above. $\hat{P}_{i}$ and $\tilde{P}_{i}$, $i=1,\ldots 7$, are the respective areas of the polygons which we are summing over. Note that $\hat{P}_{1}=\tilde{P}_{1}$ and $\hat{P}_{7}=\tilde{P}_{7}$. The process of finding all the $\hat{P}_{i}$ and $\tilde{P}_{i}$ is long but fairly elementary.

We now obtain leading order expressions for each sum in $t$. The way to do this is similar as in \cite{OreDet09}, where a more detailed explanation of the method can be found. First we substitute $t=\frac{1}{u}$, and then $n=\frac{v}{u}$ into the $\hat{P}_{i}$ and the $\tilde{P}_{i}$, such that $u$ is small and $v=\mathcal{O}\left(1\right)$. We Taylor expand $\hat{P}_{i}$ and $\tilde{P}_{i}$ into series up to order $u^2$ and then reverse the substitution by setting $v= n u$, thus effectively incorporating the large $n$ into the leading order term of each series expansion. Now each sum can be simplified into expressions involving polygamma functions of order $0$ and $1$.

The polygamma function of order $i$ is defined as the $(i + 1)$th derivative of the logarithm of the gamma function:
\begin{equation}
\Psi^{(i)}(z) = \frac{\mathrm{d}^{(i+1)}}{\mathrm{d}z^{(i+1)}} \ln\Gamma(z).
\end{equation}
The polygamma functions are of the form $z = \frac{a}{b u} +c$, where $a, b$ and $c$ are real constants, and can thus be
expanded as a Taylor series to leading order as follows:
\begin{align}
\Psi^{(0)}\left(\frac{a}{b u} + c\right) &= \ln\left(\left|\frac{a}{b u}\right|\right) + \mathcal{O}(u),\\
\Psi^{(i \geq 1)}\left(\frac{a}{b u} + c\right) &= (-1)^{(i-1)}(i-1)!\Big(\frac{b u}{a}\Big)^{i} + \mathcal{O}(u^{i+1}).
\end{align}
With these approximations at hand, we obtain expressions for the sums in (48) and (49). We only present here the first of the
approximated sums and include the rest in Appendix \textit{D}:
\begin{equation}
\sum_{n=0}^{n_{A}} \hat{P}_{1} =\sum_{n=0}^{n_{A}} \tilde{P}_{1}= \frac{(h^{+} - L)^2}{2 (2 \rho + 1) t},
\end{equation}
where we have neglected terms of order $\sim t^{-2}$.
Altogether (48) and (49) take the form:
\begin{equation}
\frac{(L-h^{+})^{2}}{4\zeta (1+\zeta)\rho t}\left[\frac{\varepsilon_{1}\rho +\varepsilon_{2}\rho^{2}+ \varepsilon_{3}\rho^{3}+\varepsilon_{4}\rho^{4}}{(2\rho+1)(2\zeta \rho-1)^{2}}+ \ln \Big((2\rho+1)^{\jmath_{1}}(2\zeta \rho -1)^{\jmath_{2}}\Big)\right],
\end{equation}
where the coefficients $\varepsilon_{i}$ ($i=1\ldots 4$) and $\jmath_{j}$ ($j=1,2$) are given in Appendix \textit{E} for both orderings $\hat{\cdot}$ and $\tilde{\cdot}$. It remains to multiply (50) by $2$ due to the horizontal symmetry of the mushroom, and normalize by $2|\partial Q|\mathcal{B}$ to obtain a probability.
The sum of expressions (39-40) and (54) depending on the value of $\zeta$, therefore gives the asymptotic contribution of the long surviving near-bouncing ball orbits $\frac{\mathcal{C}}{t}$ to the mushroom's survival probability $P(t)$.

Interestingly yet reassuringly, in the limit of $\rho\rightarrow 1$, $\zeta=2$ and the complicated expression for (54) reduces to $\frac{(L-h^{+})^{2}(3\ln3+2)}{4 t}$, which is exactly what one would expect since in this limit, the mushroom billiard is reduced to the half-stadium billiard \cite{OreDet09}. In the opposite limit where  $\rho\rightarrow 0$, $\zeta \rightarrow \infty$ the mushroom's stem shrinks and expression (54) has asymptotic expansions of
\begin{equation}
\frac{1}{t}\left(\frac{7}{2}(L-h^{+})^{2} -2(1+\zeta)(L-h^{+})^{2}\rho + 3\zeta (L-h^{+})^{2} \rho^{2}\right)+ \mathcal{O}(r^{3})
\end{equation}
and
\begin{equation}
\frac{1}{t}\left(-\frac{5}{2}(L-h^{+})^{2} +4(1+\zeta)(L-h^{+})^{2}\rho - 4(2\zeta-3) (L-h^{+})^{2} \rho^{2}\right)+ \mathcal{O}(r^{3}),
\end{equation}
for the two orderings $\hat{\cdot}$ and $\tilde{\cdot}$ respectively, indicating that the discontinuous dependence on $\zeta$ persists, hence this limit is in some sense ill-defined.

\subsection{NUMERICAL SIMULATIONS}

\begin{figure}[h]
\begin{center}
\fbox{
\includegraphics[scale=0.27]{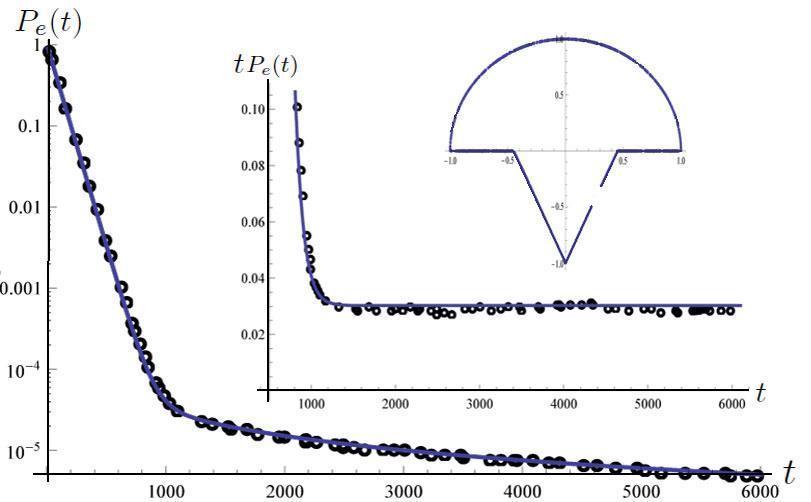}
\includegraphics[scale=0.27]{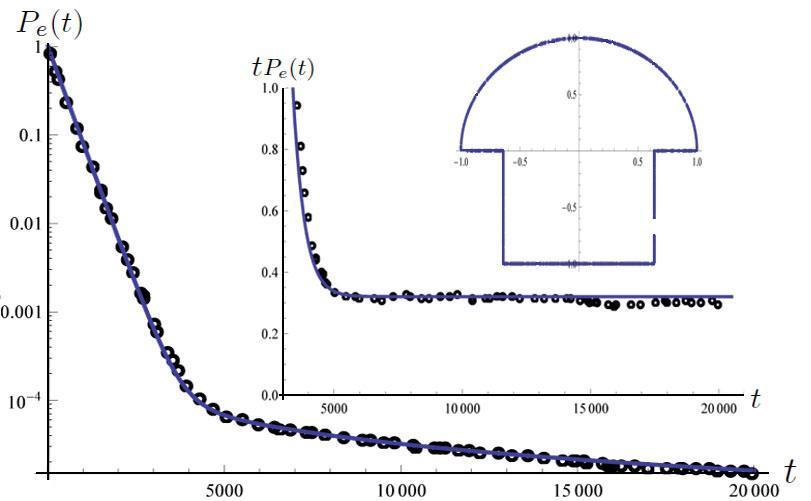}
\includegraphics[scale=0.28]{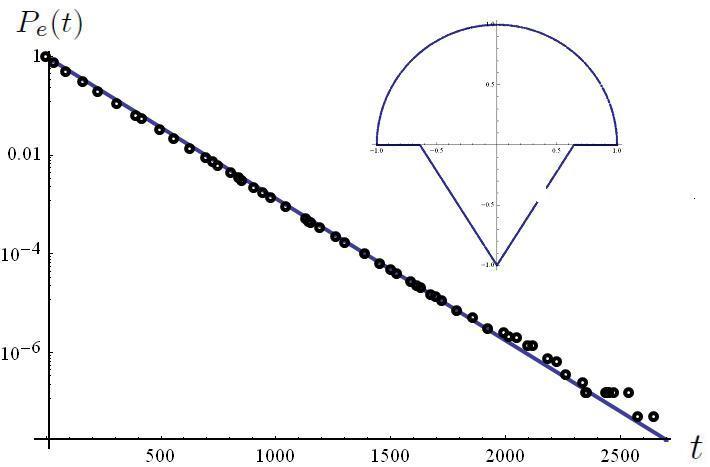}}
\caption{\label{fig:numerics} \footnotesize (Color online) Numerical simulations of $P_{e}(t)$ defined in (57) are plotted on a logarithmic scale using $10^{8}$ chaotic ICs as a function of $t$. The parameters $(r,R,L,\epsilon)$ used for the triangular stem (\emph{Left}) are $(\cos 0.3484\pi, 1, 1, 0.048)$ such that $\mathcal{S}_{\rho}=\{(20,7), (66,23), (376,131)\}$, while for the rectangular stem (\emph{Middle}) $(\cos\left(\frac{5+\sqrt{2}}{23}\pi\right), 1, 1, 0.02)$, with $h^{+}=0.3$ such that the mushroom's hat has no MUPOs. The MUPO-free mushroom (\emph{Right}) has parameters $(\cos 0.3484\pi, 1, 1, 0.0371)$. The blue curves are the analytic predictions while the numerical data correspond to the empty circles. The insets are plots of $t P_{e}(t)$ showing the agreement with the analytic expressions for $\frac{\mathcal{C}}{t}$.}
\end{center}
\end{figure}

Having obtained exact leading order analytic expressions for all the parameters appearing in (25) we now numerically test their validity by plotting the conditional probability $P_{e}(t)$ that a particle survives up to time $t$ given that the particle is chosen uniformly from the ergodic component of the billiard flow (see Figure ~\ref{fig:numerics}).
\begin{equation}
P_{e}(t)=(P(t)-\mathcal{A})/\mathcal{B}= e^{-\bar{\gamma}t}+\frac{\mathcal{C}}{t} + \mathcal{O}\left(\frac{1}{t^{2}}\right).
\end{equation}
The plots are purposely chosen (from many more) to portray and verify the results obtained in the present paper. Three different mushrooms are simulated: one with MUPOs present only in the hat (\emph{Left}), one with bouncing ball orbits in the stem and a MUPO-free hat (\emph{Middle}), and one with no MUPOs at all (\emph{Right}). Different hole sizes give different exponential escape rates $\bar{\gamma}$ and in turn cross-over times to a power law decay. The empty black circles in the plots correspond to the numerical data while the blue curves give the analytic predictions of (57). Although each simulation consists of $10^{8}$ chaotic ICs, we were not able to pick up any power law decay in the MUPO-free mushroom (\emph{Right}) and hence any clues of stickiness.

\section{CONCLUSIONS}

In this paper, we have attempted to quantify the stickiness observed in the mushroom billiard by placing a hole in its ergodic component and looking at the survival probability function $P(t)$ at long times (see eq (25)). Our analytic predictions are in good agreement with the numerical simulations performed and therefore confirm that $P(t)\sim \frac{\mathcal{C}}{t}$ for long enough times. Also, their good agreement with the constants $\mathcal{C}$ derived in sections III A. and B. for MUPOs present in the hat and in the stem respectively, implies that these MUPOs are indeed the primary causes of the power-law decay. This observation in turn applies to the Poincar\'{e} recurrence times distribution $\mathcal{Q}(t)$ studied in Ref \cite{Altman09} and the rate of mixing of the ergodic component \cite{Chernov08}.

The explicit expressions obtained here for $\mathcal{C}$, allow one not only to predict but also to calibrate the asymptotic behaviour of $P(t)$. Also, we have shown that these distributions as well as the overall existence of MUPOs in the hat are sensitive to the system's control parameter $\rho=\frac{r}{R}$, whilst only the near-bouncing ball orbits' contribution to $P(t)$ depends on the hole's position and size. The reason for this is that the hole intersects the sticky region generated in phase space by the period-2 bouncing ball orbits. This creates a fictitious, time dependent `island of stability' in the mushroom's ergodic component. Although orbits in it are unstable, they only experience up to one non-linear collision process before escaping, thus allowing us to approximate their occupancy in phase space with polygonal `spikes' which we could then integrate over. In the case of the MUPOs in the mushroom's hat, we could easily bound the long surviving orbits by assuming that they will escape exponentially fast once in the stem.
It is expected that the methods used here can be further generalised and applied to other mushrooms with elliptical hats for instance, or even to other billiards such as the annular or drive-belt stadium billiards where polygonal type MUPOs act as scaffolding for sticky orbits to cling onto.

\begin{figure}[h]
\begin{center}
\fbox{
\includegraphics[scale=0.25]{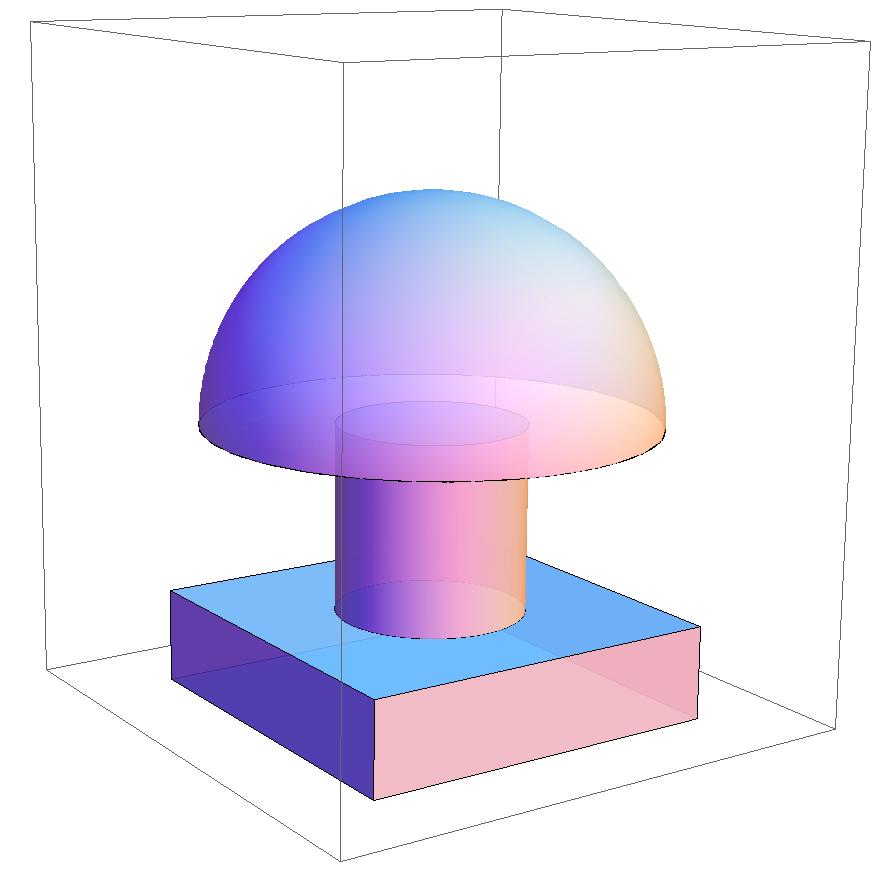}}
\caption{\label{fig:mush3d} \footnotesize (Color online) A three-dimensional mushroom billiard.}
\end{center}
\end{figure}

A major result of this paper is the introduction of a zero measure set which describes MUPO-free mushrooms (see section II B.) which to the best of our knowledge posses the simplest mixed phase space in two dimensions. The interesting connection between mushrooms and the diophantine approximation discovered by Altmann \cite{AltmanPhD} cannot be directly exploited due to a sensitive parity dependence of the periodic orbits of the mushroom. We have overcome this complication by considering a generalized mushroom with a variable sized hat and triangular stem (see Figure ~\ref{fig:mushroom2}). This allowed us to efficiently use properties of continued fractions and characterize a subset of the infinitely many MUPO-free mushrooms. We thus obtained upper bounds for MUPO-free and finitely sticky irrational mushrooms and also gave an explicit example of a MUPO-free mushroom billiard (see section II C.).

Unlike the non-sticky elliptical mushrooms mentioned in the introduction, the MUPO-free mushroom with a circular type hat exhibits a reduced amount of stickiness (larger scaling exponent). This is due to the difficulty in `finding' the foot of the mushroom by orbits which are just inside the dashed circular arc of radius $r$. Hence, the MUPO-free mushroom, although not generic, is mathematically important as it should display a polynomial decay of correlations of order $\sim t^{-2}$.

The results of section II also apply in the case of a three dimensional mushroom billiard with a hemispherical hat of radius $R$, a cylindrical stem of radius $r$ and height $h>0$ and a cuboidal pedestal of base length $l\geq2r$ to break angular momentum conservation (see Figure ~\ref{fig:mush3d}). This is because orbits inside a three-dimensional spherical billiard always lie in the same two-dimensional plane which also contains the center of the corresponding sphere.
The remaining (zero-measure) MUPOs in such a system are of the bouncing ball type and are found both in stem and pedestal. We conjecture that the corresponding mushroom has a sharply divided phase space with a MUPO-free hat \cite{Bun10}.

Finally, one would expect to see the classical dynamical features caused by the removal of MUPOs to appear in the analogous quantum system in accordance with Bohr's correspondence principle. Two obvious candidates for such an observation are the localization of wave-functions (scars) \cite{Barnett07} and the effect on the dynamical tunneling rates \cite{Backer08}, thus providing a new model for quantum chaos \cite{Stockmann99} as well as experimentalists with applications which range from semiconductor nano-structures to dielectric micro-cavities.



\section*{ACKNOWLEDGEMENTS}
We would like to thank Sandro Bettin, Leonid Bunimovich and Alan Haynes for helpful discussions and OG's EPSRC Doctoral Training Account number SB1715.

\appendix
\section{A Bound on $R_{2}$}
To obtain a bound on $R_{2}$ in eq (8) we shall use the remainder term from Taylor's theorem several times. Taylor's theorem states that if $f$ is a function which is $n$ times differentiable on the closed interval $[a, x]$ and $n + 1$ times differentiable on the open interval $(a, x)$, then
\est{f(x)=f(a)+\frac{f'(a)}{1!}(x-a)+\frac{f''(a)}{2!}(x-a)^{2}+\ldots +\frac{f^{(n)}(a)}{n!}(x-a)^{n}+R_{n}(x),
}
where $R_{n}(x)=\frac{f^{(n+1)}(\xi)}{(n+1)!}(x-a)^{n+1}$ for some $a<\xi<x$.

Let $c=\cos\pr{\alpha \pi\frac pq}\leq\cos\pi\vartheta^{*}$, and $\varepsilon=c\pr{\frac1{\cos\pr{\frac{\alpha \pi}q}}-1}$ such that eq (6) corresponds to
\est{c\leq\rho<c+\varepsilon.}
We have to compute
\est{
\arccos\pr{c+\varepsilon}=\arccos\pr{\frac{\cos\pr{\alpha \pi\frac pq}}{\cos\pr{\frac{\alpha \pi}q}}}.
}
The Taylor expansion of $\arccos\pr{c+\varepsilon}$ at $\varepsilon=0$ is
\es{\label{1}
\arccos\pr{c+\varepsilon}=\arccos(c)-\frac\varepsilon{\sqrt{1-c^2}}+A_1
}
with $A_1=-\frac{c+x}{2\pr{1-\pr{c+x}^2}^{\frac32}}\varepsilon^2<0$ and $0\leq x \leq \varepsilon$.
Now, for $q\rightarrow\infty$, we have
\es{\label{2}
\frac\varepsilon c=\pr{\frac1{\cos(\frac{\alpha \pi}q)}-1}&=\frac{\alpha^2\pi^2}{2q^2}+B_2=B_0
}
with $q^3 B_2=\frac{\sin(\alpha \pi z)^3\alpha ^3\pi^3}{\cos(\alpha \pi z)^4}+\frac{5\sin(\alpha \pi z)\alpha ^3 \pi^3}{6\cos(\alpha \pi z)^2}$ and $q B_0=\frac{\sin(\alpha \pi w)\alpha\pi}{\cos(\alpha\pi w)^2}$ for some $0\leq z,w\leq \frac1q$. Since for $q\geq Q$ we have $1<\frac1{\cos(\alpha \pi z)}\leq \frac1{\cos\pr{\frac{\alpha \pi }{Q}}}$ and $0< \sin(\alpha \pi z)\leq \sin\pr{\frac{\alpha \pi }{Q}}.$ Hence
\est{
0<q^3 B_2\leq \pr{\frac{\sin(\frac{\alpha \pi}{Q})^2}{\cos(\frac{\alpha \pi}{Q})^4}+\frac{5}{6\cos(\frac{\alpha \pi}{Q})^2}}\alpha^3\pi^3\sin\pr{\frac{\alpha \pi}{q}}\leq \pr{\frac{\sin(\frac{\alpha \pi}{Q})^2}{\cos(\frac{\alpha \pi}{Q})^4}+\frac{5}{6\cos(\frac{\alpha \pi}{Q})^2}}\frac{\alpha^4\pi^4}{q},
}
and similarly
\est{
0<q B_0\leq \frac{\alpha^2\pi^2}{\cos(\frac{\alpha \pi}{Q})^{2}{q}}.
}
This, together with~\eqref{1} and~\eqref{2}, gives
\est{
\arccos(c+\varepsilon)=\arccos(c)-\frac{\alpha^2c\pi^2}{2\sqrt{1-c^2}}\frac{1}{q^2}+C_1+A_1
}
where $C_1=-\frac{c B_2}{\sqrt{1-c^2}}<0$ is bounded by
\est{
|C_1|\leq \pr{\frac{\sin(\frac{\alpha \pi}{Q})^2}{\cos(\frac{\alpha \pi}{Q})^4}+\frac{5}{6\cos(\frac{\alpha \pi}{Q})^2}}\frac{\alpha^4\pi^4}{q^{4}}\frac{c}{\sqrt{1-c^2}} \leq \pr{\frac{\sin(\frac{\alpha \pi}{Q})^2}{\cos(\frac{\alpha \pi}{Q})^4}+\frac{5}{6\cos(\frac{\alpha \pi}{Q})^2}}\frac{\alpha^4\pi^4}{q^{4}}\frac{\cos\pi\vartheta^{*}}{\sqrt{1-(\cos\pi\vartheta^{*})^2}}
}
and $A_1 <0$ by
\est{
|A_1|&\leq\frac{c(1+B_0)}{2\pr{1-c^2(1+B_0)^2}^\frac32}c^2B_0^2 \leq \pr{1+\frac{\alpha^{2}\pi^{2}}{\cos(\frac{\alpha\pi}{Q})^{2}Q^{2}}}\frac{\alpha^4\pi^4c^3}{2 \cos(\frac{\alpha \pi}{Q})^{4}\pr{1-c^2(1+B_0)^2}^\frac{3}{2} q^4} \\
&\leq	 \pr{1+\frac{\alpha^{2}\pi^{2}}{\cos(\frac{\alpha\pi}{Q})^{2}Q^{2}}}\frac{\alpha^4\pi^4(\cos\pi\vartheta^{*})^3}{2 \cos(\frac{\alpha \pi}{Q})^{4}\pr{1-(\cos\pi\vartheta^{*})^2(1+\frac{\alpha^2\pi^2}{\cos(\frac{\alpha \pi}{Q})^{2}{q^{2}}})^2}^\frac{3}{2} q^4}
}
for $q>\frac{\alpha\pi}{\cos(\frac{\alpha \pi}{Q})}\sqrt{\frac{\cos\pi\vartheta^{*}}{(1-\cos\pi\vartheta^{*})}}$. In the same way
\es{\label{3}
\frac{\alpha^2c\pi^2}{2\sqrt{1-c^2}}\frac{1}{q^2}=\frac{\alpha^2 \pi^2 \cos\pi\vartheta^{*}}{2\sqrt{1-c^2}}\frac{1}{q^2}+C_2,
}
with $C_{2}<0$ and bounded by
\est{
|C_2|&\leq\frac{\alpha^2 \pi^2 \varepsilon }{2\sqrt{1-c^2}}\frac{1}{q^2} \leq \frac{\alpha^2 \pi^2 B_0 \cos\pi\vartheta^{*} }{2\sqrt{1-(\cos\pi\vartheta^{*})^2}}\frac{1}{q^2}\leq\frac{\alpha^4\pi^4\cos\pi\vartheta^{*}}{2\cos(\frac{\alpha \pi}{Q})^2\sqrt{1-(\cos\pi\vartheta^{*})^2}}\frac{1}{q^4}.
}
Finally we must bound
\est{
\frac1{\sqrt{1-c^2}}=\frac1{\sqrt{1-(\cos\pi\vartheta^{*})^2 + ((\cos\pi\vartheta^{*})^2-c^{2})}}=\frac1{\sqrt{1-(\cos\pi\vartheta^{*})^2 + \nu}},
}
where we wrote $\nu=(\cos\pi\vartheta^{*})^{2}-c^2$ so that we may expand for $\nu$ small
\est{\frac1{\sqrt{1-(\cos\pi\vartheta^{*})^2 + \nu}}= \frac1{\sqrt{1-(\cos\pi\vartheta^{*})^2}}+D_0
}
where $D_0=-\frac{\nu}{2(1-(\cos\pi\vartheta^{*})^2+y)^{\frac32}}$, with $0\leq y \leq\nu$ and so it is bounded by
\est{
|D_0|\leq\frac{\nu}{2\kappa^{\frac32}}\leq\frac{2\varepsilon\cos\pi\vartheta^{*}}{2\pr{1-(\cos\pi\vartheta^{*})^2}^\frac32}\leq\frac{\alpha^2\pi^2(\cos\pi\vartheta^{*})^2}{\cos(\frac{\alpha \pi}{Q})^2\pr{1-(\cos\pi\vartheta^{*})^2}^\frac32}\frac1{q^2}.
}
Therefore for ~\eqref{3} we have
\est{
\frac{\alpha^2 \pi^2 c}{2\sqrt{1-c^2}}\frac{1}{q^2} = \frac{\alpha^2 \pi^2 \cos\pi\vartheta^{*}}{2\sqrt{1-(\cos\pi\vartheta^{*})^2}}\frac{1}{q^2}+C_2+C_3,
}
where
\est{
|C_3|\leq\frac{\alpha^4\pi^4(\cos\pi\vartheta^{*})^3}{2\cos(\frac{\alpha \pi}{Q})^2\pr{1-(\cos\pi\vartheta^{*})^2}^\frac32}\frac1{q^4}.
}
Putting everything together, for $q\geq \mathrm{max}\pr{Q,\frac{\alpha\pi}{\cos(\frac{\alpha \pi}{Q})}\sqrt{\frac{\cos\pi\vartheta^{*}}{(1-\cos\pi\vartheta^{*})}}}$ we have
\est{
\arccos(c+\varepsilon)&=\arccos(c)-\frac{\alpha^2 c\pi^2}{2\sqrt{1-c^2}}\frac{1}{q^2}+C_1+A_1\\
&=\arccos(c)-\frac{\alpha^2 \pi^2 \cot\pi\vartheta^{*}}{2 q^2}+C_1+A_1+C_2+C_3,
}
where the remainders have magnitudes bounded by
\est{
|C_1|&\leq \pr{\frac{\sin(\frac{\alpha \pi}{Q})^2}{\cos(\frac{\alpha \pi}{Q})^4}+\frac{5}{6\cos(\frac{\alpha \pi}{Q})^2}} \frac{\alpha^4\pi^4\cot\pi\vartheta^{*}}{q^{4}},\\
|A_1|&\leq\pr{1+\frac{\alpha^{2}\pi^{2}}{\cos(\frac{\alpha\pi}{Q})^{2}Q^{2}}}\frac{\alpha^4\pi^4(\cos\pi\vartheta^{*})^3}{2 \cos(\frac{\alpha \pi}{Q})^{4}\pr{1-(\cos\pi\vartheta^{*})^2(1+\frac{\alpha^2\pi^2}{\cos(\frac{\alpha \pi}{Q})^{2}{q^{2}}})^2}^\frac{3}{2} q^4},\\
|C_2|&\leq\frac{\alpha^4\pi^4\cot\pi\vartheta^{*}}{2\cos(\frac{\alpha \pi}{Q})^2 q^4},\\
|C_3|&\leq\frac{\alpha^4\pi^4(\cot\pi\vartheta^{*})^3}{2\cos(\frac{\alpha \pi}{Q})^2 q^4},
}
and are all negative.

\section{Corners of the polygonal Area}
The corners of the polygons (for fixed $n$) as shown in Figure ~\ref{fig:hyperbolas} are found by solving for the intersections of the various curves and lines obtained from equations (45) and (46):
\begin{align}
A_{x_{i}}&=\frac{-2 (h^{+}\rho (1 + 2\rho) (1 + n +\rho + 2 n\rho) - L (1 +\rho) (1 + 2 (2 + n)\rho + (2 + 4 n)\rho^2))+ R (-2 L \rho R - L R) t}{2 + 2 (4 + n)\rho + (6 + 4 n)\rho^2 + R (-2 \rho R - R) t},\\
A_{\theta_{i}}&=\frac{(h^{+}-L) (1 + 2 \rho)^2}{ 2 R (1 + \rho (4 + n + 3 \rho + 2 n \rho))-t (1 + 2 \rho) },\\
B_{x_{i}}&=\frac{-2 h^{+} R \rho (1 + 2 \rho) (1 + n + \rho + 2 n \rho) +  L (-2 (1 + \zeta) R + t + 4 (1 + 3 \zeta + 2 (1 + \zeta) n) R \rho^3)}{(2 \zeta \rho -1) (-t (1 + 2 \rho) + 2 R (1 + \zeta + 2 \zeta \rho + \rho (4 + n + 3 \rho + 2 n \rho)))}
\nonumber \\
&\quad+ \frac{2 \rho (t - \zeta t (1 + 2 L \rho) + R (-3 + 2 \zeta^2 + 2 L (n + \zeta (4 + 2 \zeta + n)) \rho))}{(2 \zeta \rho -1) (-t (1 + 2 \rho) +
2 R (1 + \zeta + 2 \zeta \rho + \rho (4 + n + 3 \rho + 2 n \rho)))},\\
B_{\theta_{i}}&=\frac{(h^{+} - L) (2 \rho + 1)^2}{( 2 \zeta \rho-1) (-t (1 + 2 \rho) +
2 R (1 + \zeta + 2 \zeta \rho + \rho (4 + n + 3 \rho + 2 n \rho)))},\\
C_{x_{i}}&=L - \frac{2 (h^{+} - L) R \rho (2 \zeta \rho-1) (\zeta \rho + n (2 \zeta \rho-1))}{(1 + 2 \rho) (t - 2 \zeta t \rho + 2 R (-1 - \zeta + (-1 + 2 \zeta (1 + \zeta) - n) \rho + \zeta (3 + 2 n) \rho^2))},\\
C_{\theta_{i}}&=\frac{(h^{+} - L) (1 - 2 \zeta \rho)^2}{(1 + 2 \rho) (t - 2 \zeta t \rho +2 R (-1 - \zeta + (-1 + 2 \zeta (1 + \zeta) - n) \rho +
\zeta (3 + 2 n) \rho^2))},\\
D_{x_{i}}&=\frac{2 (L + h^{+} n) R \rho -L t + 4 \zeta^2 R \rho ((h^{+} - L) (1 + 2 n) \rho^2 -L) + 2 \zeta (L R + L t \rho - (h^{+} - 2 L (n-1) + 4 h^{+} n) R \rho^2)}{-t - 4 \zeta^2 R \rho + 2 (1 + n) R \rho + 2 \zeta (R + t \rho - (3 + 2 n) R \rho^2)},\\
D_{\theta_{i}}&=\frac{(h^{+} - L) (1 - 2 \zeta \rho)^2}{-2 \zeta R + t - 2 ((1 - 2 \zeta^2 + n) R + \zeta t) \rho +
 2 \zeta (3 + 2 n) R \rho^2},\\
E_{x_{i}}&= \frac{L t + 4 \zeta^2 L R \rho - 2 (h^{+} + L + h^{+} n) R \rho + 2 \zeta L ( R (-1 + 2 (2 + n) \rho^2)-t \rho)}{(2 \zeta \rho -1) ( 2 R (\zeta + (2 + n) \rho)-t)},\\
E_{\theta_{i}}&= \frac{(h^{+} - L)}{(-1 + 2 \zeta \rho) (-t + 2 R (\zeta + (2 + n) \rho))},\\
F_{x_{i}}&= \frac{L (2 R (1 + \rho) (1 + 2 (1 + n) \rho) - t (1 + 2 \rho)) - 2 h^{+} n R \rho}{(1 + 2 \rho) (2 R (1 + \rho + n \rho) - t)},\\
F_{\theta_{i}}&= \frac{(h^{+} - L) }{(1 + 2 \rho) (-t + 2 R (1 + \rho + n \rho))},\\
G_{x_{i}}&= L,\\
G_{\theta_{i}}&= 0,
\end{align}
where $\zeta=\left\lceil\frac{R}{r}\right\rceil$ and $\rho=\frac{r}{R}$.

\section{Values of $n$ when corners hit the hole}
The upper and lower limits of the sums in expressions (48-49) are the solutions for $n$ when the $x_{i}$ coordinate of the corners $A-F$ exceeds $h^{+}$:
\begin{align}
n_{B}&= \left\lfloor\frac{(2+\zeta-\frac{t}{R})+\left(6-4\zeta^{2}+ \frac{t}{R}(2\zeta-2)\right)\rho   +(-16\zeta-8\zeta^{2}+4\zeta\frac{t}{R})\rho^{2}+(-4-12\zeta)\rho^{3}}{4(1+\zeta)(2\rho+1)\rho^{2}}\right\rfloor ,\\
n_{C}&= \left\lfloor\frac{(2+\zeta-\frac{t}{R})+\left(6-4\zeta^{2}+ \frac{t}{R}(2\zeta-2)\right)\rho   +(4-12\zeta-8\zeta^{2}+4\zeta\frac{t}{R})\rho^{2}+(-12\zeta-4\zeta^{2})\rho^{3}}{4(1+\zeta)(2\zeta\rho-1)\rho^{2}}\right\rfloor ,\\
n_{D}&= \left\lfloor\frac{2\zeta-\frac{t}{R}+(2-4\zeta^{2}+2\zeta \frac{t}{R})\rho - 4\zeta \rho^{2}-4\zeta^{2}\rho^{3}}{4\zeta (2\zeta\rho-1)\rho^{2}}\right\rfloor ,\\
n_{E}&= \left\lfloor\frac{2\zeta-\frac{t}{R}+(2-4\zeta^{2}+2\zeta\frac{t}{R})\rho-8\zeta\rho^{2}}{4\zeta\rho^{2}}\right\rfloor ,\\
n_{F}&= \left\lfloor\frac{-2+\frac{t}{R}+(-6+\frac{t}{R})\rho-4\rho^{2}}{4\rho(\rho+1)}\right\rfloor,
\end{align}
where $\zeta=\left\lceil\frac{R}{r}\right\rceil$ and $\rho=\frac{r}{R}$.

\section{Leading order approximations of Sums in eq (48-49)}
\begin{align}
\sum_{n=n_{A}+1}^{n_{B}} \hat{P}_{2} &= \frac{(h^{+} - L)^2}{4 \rho t}\left[\frac{(4 + 2 \zeta ) \rho + (-8 \zeta - 2 \zeta^2 ) \rho^2 + 4 \zeta^2 \rho^3}{(2 \zeta \rho - 1) (2 \rho + 1)} + \ln \left(2 \zeta \rho - 1\right)   \right],\\
\sum_{n=n_{B}+1}^{n_{D}} \hat{P}_{3} &= \frac{(h^{+} - L)^2}{4 (1 + \zeta) t \rho}\left[\frac{2 \rho (-1 + 2 \zeta ( \zeta \rho)-1) (2 \rho + \zeta (-1 + \rho + 2 \zeta \rho - 2 \rho^2))}{(1 + 2 \rho) (1 - 2 \zeta \rho)^2} + (2 + \zeta) \ln \left(\frac{2 \rho + 1}{(2 \zeta \rho - 1)^2}\right)   \right],\\
\sum_{n=n_{D}+1}^{n_{E}} \hat{P}_{4} &= \frac{(h^{+} - L)^2}{4 \zeta (1 + \zeta) t \rho}\left[\frac{-2 \rho (1 + \zeta - \zeta^2 + 2 \zeta^3 \rho) (1 + \rho (1 + 2 \zeta (-1 + \zeta \rho)))}{(1 + 2 \rho) (1 - 2 \zeta \rho)^2} + (1 + 3 \zeta + \zeta^2) \ln \left(2  \rho - 1\right)\right],\\
\sum_{n=n_{E}+1}^{n_{C}} \hat{P}_{5} &= \frac{(h^{+} - L)^2}{4 (1 + \zeta) t \rho}\left[\frac{2 \rho (-1 + 2 \zeta (-1 + \zeta \rho)) (2 \rho + \zeta (-1 + \rho + 2 \zeta \rho - 2 \rho^2))}{(1 + 2 \rho) (1 - 2 \zeta \rho)^2} + (2 + \zeta) \ln \left(\frac{2 \rho + 1}{(2 \zeta \rho - 1)^2}\right)   \right],\\
\sum_{n=n_{C}+1}^{n_{F}} \hat{P}_{6} &= \frac{(h^{+} - L)^2}{4 \rho t}\left[\frac{2 \rho (\zeta \rho - 1) (\zeta (2 \rho - 1) - 2)}{(1 + 2 \rho) (2 \zeta \rho - 1)} + \ln \left(2 \zeta \rho - 1\right)   \right],\\
\sum_{n=n_{F}+1}^{\infty} \hat{P}_{7} &= \frac{(h^{+} - L)^2 (\rho + 1) }{(2 \rho + 1) t},\\
\sum_{n=n_{A}+1}^{n_{D}} \tilde{P}_{2} &= \frac{(h^{+} - L)^2}{4 \rho t}\left[\frac{(2 - 2 \zeta ) \rho + (2 - 4 \zeta + 2 \zeta^2) \rho^2}{(2 \zeta \rho - 1) (2 \rho + 1)} + \ln \left(\frac{2 \rho + 1}{2 \zeta \rho - 1}\right)   \right],\\
\sum_{n=n_{D}+1}^{n_{B}} \tilde{P}_{3} &= \frac{(h^{+} - L)^2}{4 \zeta t \rho}\left[\frac{ 2 \rho (1 + \rho + 2 \zeta^2 \rho - \zeta (1 + \rho)) (1 - 2 \zeta (\zeta \rho - 1))}{(1 + 2 \rho) (1 - 2 \zeta \rho)^2} + (1 + \zeta) \ln \left( \frac{(2 \zeta \rho - 1)^2}{(2 \rho + 1)}\right)\right],\\
\sum_{n=n_{B}+1}^{n_{C}} \tilde{P}_{4} &= \frac{(h^{+} - L)^2}{4 \zeta (1 + \zeta) t \rho}\left[\frac{ 4 (1 + \zeta) \rho (1 - (\zeta - 1) \rho) (1 + \zeta - \zeta^2 + 2 \zeta^3 \rho)}{(1 + 2 \rho) (1 - 2 \zeta \rho)^2} + (1 + 3 \zeta + \zeta^2) \ln \left(\frac{(2 \zeta \rho - 1)^2}{(2 \rho + 1)^2}\right)\right],\\
\sum_{n=n_{C}+1}^{n_{E}} \tilde{P}_{5} &= \frac{(h^{+} - L)^2}{4 \zeta t \rho}\left[\frac{ 2 \rho (1 + \rho + 2 \zeta^2 \rho - \zeta (1 + \rho)) (1 - 2 \zeta (\zeta \rho - 1))}{(1 + 2 \rho) (1 - 2 \zeta \rho)^2} + (1 + \zeta) \ln \left(\frac{(2 \zeta \rho - 1)^2}{2 \rho + 1}\right)   \right],\\
\sum_{n=n_{E}+1}^{n_{F}} \tilde{P}_{6} &= \frac{(h^{+} - L)^2}{4 \rho t}\left[\frac{(2 - 2 \zeta) \rho + (2 - 4 \zeta + 2 \zeta^2) \rho^2}{(1 + 2 \rho) (2 \zeta \rho - 1)} + \ln \left(\frac{2\rho+1}{2 \zeta \rho - 1}\right)   \right],\\
\sum_{n=n_{F}+1}^{\infty} \tilde{P}_{7} &= \frac{(h^{+} - L)^2 (\rho + 1) }{(2 \rho + 1) t},
\end{align}
where $\zeta=\left\lceil\frac{R}{r}\right\rceil$ and $\rho=\frac{r}{R}$.

\section{Coefficients of eq (54)}
The coefficients of equation (54) for the two different orderings $\hat{\cdot}$ and $\tilde{\cdot}$ of $n_{A}-n_{F}$ are
\begin{align}
&\hat{\varepsilon}_{1}= -2 - 4 \zeta + 4 \zeta^3  \hspace{52.5 mm} \tilde{\varepsilon}_{1}= 8 + 18 \zeta + 2 \zeta^2 - 8 \zeta^3  \\
&\hat{\varepsilon}_{2}= -2 (1 + \zeta + \zeta^2 + 2 \zeta^3 + 6 \zeta^4 ) \hspace{32.7 mm} \tilde{\varepsilon}_{2}= 8 + 12 \zeta - 20 \zeta^2 + 24 \zeta^4 \\
&\hat{\varepsilon}_{3}= -12 \zeta^2 - 8 \zeta^3 - 4 \zeta^4 + 8 \zeta^5 \hspace{36.3 mm} \tilde{\varepsilon}_{3}= -16 \zeta^2 + 8 \zeta^3 + 8 \zeta^4 - 16 \zeta^5\\
&\hat{\varepsilon}_{4}= 16 \zeta^3 + 16 \zeta^4 + 8 \zeta^5 \hspace{46.8 mm} \tilde{\varepsilon}_{4}= 16 (\zeta^3 + \zeta^4)\\
&\hat{\jmath}_{1}= 1 + 7 \zeta + 3 \zeta^2 \hspace{55.8 mm} \tilde{\jmath}_{1}= -4 - 8 \zeta - 2 \zeta^2\\
&\hat{\jmath}_{2}= -6 \zeta - 2 \zeta^2 \hspace{59.25 mm} \tilde{\jmath}_{2}= 6 + 12 \zeta + 4 \zeta^2
\end{align}
where $\zeta=\left\lceil\frac{R}{r}\right\rceil$.

\end{document}